\documentclass[12pt,a4paper]{amsart}

\usepackage{color}
\usepackage{amscd,amssymb,amsopn,amsmath,amsthm,graphics,amsfonts,enumerate,verbatim,calc}
\usepackage[dvips]{graphicx}

\newtheorem{Theorem}{Theorem}[section]
\newtheorem{Proposition}[Theorem]{Proposition}
\newtheorem{Lemma}[Theorem]{Lemma}

\newtheorem{Remark}[Theorem]{Remark}
\newtheorem{Example}[Theorem]{Example}
\newtheorem{Definition}[Theorem]{Definition}

\def\RMN#1{\uppercase\expandafter{\romannumeral#1}}

\newcommand{\spec}{\mathop{\rm Spec}\nolimits}
\newcommand{\proj}{\mathop{\rm Proj}\nolimits}
\newcommand{\depth}{\mathop{\rm depth}\nolimits}

\newcommand{\assh}{\mathop{\rm Assh}\nolimits}

\newcommand{\height}{\mathop{\rm ht}\nolimits}

\newcommand{\grade}{\mathop{\rm grade}\nolimits}

\newcommand{\gp}{\mathfrak{p}}

\newcommand{\gm}{\mathfrak{m}}

\newcommand{\ga}{\mathfrak{a}}
\newcommand{\gb}{\mathfrak{b}}
\newcommand{\gc}{\mathfrak{c}}

\newcommand{\zz}{\mathbb{Z}}

\newcommand{\qq}{\mathbb{Q}}
\newcommand{\fp}{{\mathbb{F}_p}}

\newcommand{\rees}[1]{\mathcal{R}(#1)}
\newcommand{\gr}[1]{\mathcal{G}(#1)}

\newcommand{\srees}[1]{\mathcal{R}_s(#1)}
\newcommand{\serees}[1]{\mathcal{R}'_s(#1)}
\newcommand{\sgr}[1]{\mathcal{G}_s(#1)}
\newcommand{\ass}{\mathop{\rm Ass}\nolimits}
\newcommand{\hc}[3]{\mathrm{HC}(#1; #2, #3)}

\newcommand{\mult}[2]{{\mathrm e}_{#1}(#2)} %multiplicity
\newcommand{\length}[2]{\ell_{#1}(\,#2\,)} %length

\newcommand{\ra}{\longrightarrow}

\begin{document}
\title[Infinitely generated symbolic Rees rings]{
Infinitely generated symbolic Rees rings of
space monomial curves having negative curves }
\author{Kazuhiko Kurano and Koji Nishida}
\date{}
\maketitle

\begin{abstract}
In this paper, we shall study finite generation
of symbolic Rees rings of the defining ideal ${\frak p}$ of
the space monomial curve $(t^a, t^b, t^c)$
for pairwise coprime integers $a$, $b$, $c$. 
Suppose that the base field is of characteristic $0$ and
the above ideal ${\frak p}$ is minimally generated by
three polynomials.
In Theorem~\ref{Kmain1}, under the assumption that the homogeneous element $\xi$ of the minimal degree 
in  ${\frak p}$ is the negative curve,
we determine the minimal degree of an element $\eta$ such that
the pair $\{ \xi, \eta \}$ satisfies Huneke's criterion 
in the case where the symbolic Rees ring is Noetherian.
By this result, we can decide whether the symbolic Rees ring 
$\srees{{\frak p}}$ is Notherian using computers.
We give a necessary and sufficient conditions for finite generation of 
the symbolic Rees ring of ${\frak p}$ in Proposition~\ref{Kmain2}
under some assumptions.
We give an example of an infinitely generated symbolic Rees ring of ${\frak p}$
in which the homogeneous element of the minimal degree 
in  ${\frak p}^{(2)}$ is the negative curve in Example~\ref{N5f}.
We give a simple proof to (generalized) Huneke's criterion.
\end{abstract}

\section{Introduction}

Let ${\frak p}_K(a, b, c)$ be the defining ideal of the space monomial
curve $(t^a, t^b, t^c)$ in $K^3$, where $K$ is a field.
The ideal ${\frak p}_K(a, b, c)$ is generated by at most three binomials in $K[x,y,z]$
(Herzog~\cite{Her}).
The symbolic Rees rings of space monomial primes  are deeply studied by many authors. 
Huneke~\cite{Hu} and Cutkosky~\cite{Cu} developed criterions 
for finite generation of such rings. 
In 1994,
Goto-Nishida-Watanabe~\cite{GNS} first found examples of infinitely generated
symbolic Rees rings of space monomial primes.
Recently, using toric geometry,
Gonz\'ales-Karu~\cite{GK} found some sufficient conditions for 
the symbolic Rees rings of space monomial primes to be infinitely generated.

Cutkosky~\cite{Cu} found the geometric meaning of the symbolic Rees rings of space monomial primes. 
Let ${\Bbb P}_K(a,b,c)$ be the weighted projective surface with degree $a$, $b$, $c$.
Let $X_K(a,b,c)$ be the blow-up at a point in the open orbit of the toric variety ${\Bbb P}_K(a,b,c)$.
Then, the Cox ring of $X_K(a,b,c)$ is isomorphic to the extended symbolic Rees ring
of the space monomial prime ${\frak p}_K(a, b, c)$.
Therefore, the symbolic Rees ring of the space monomial prime ${\frak p}_K(a, b, c)$
is finitely generated if and only if the Cox ring of $X_K(a,b,c)$ is finitely generated,
that is, $X_K(a,b,c)$ is a Mori dream space.
A curve $C$ on $X_K(a,b,c)$ is called the negative curve if $C^2 < 0$ and
$C$ is different from the exceptional curve $E$. 
Here suppose $\sqrt{abc} \not\in {\Bbb Q}$.
Cutkosky~\cite{Cu} proved that the symbolic Rees ring of the space monomial prime ${\frak p}_K(a, b, c)$ is finitely generated if and only if
the following two conditions are satisfied:
\begin{itemize}
\item[(1)]
There exists a negative curve $C$.
\item[(2)]
There exists a curve $D$ on $X_K(a,b,c)$ such that $C \cap D =\emptyset$.
\end{itemize}

All known examples (\cite{GNW}, \cite{GK}) of the infinitely generated symbolic Rees rings of  ${\frak p}_K(a, b, c)$ satisfy
the following conditions:
\begin{itemize}
\item[(I)]
there exists a negative curve $C$ such that $C.E = 1$.
\item[(II)]
the characteristic of $K$ is $0$
\end{itemize}

In this paper, we give an example of an infinitely generated symbolic Rees ring such that there exists the negative curve $C$
with $C.E = 2$.
Furthermore, in the case where both (I) and (II) as above are satisfied,
we determine the minimal value of the degree of the  curve $D$ 
which satisfies the condition~(2) as above.

The existence of negative curves is a very difficult problem,
that is deeply related to Nagata conjecture (Proposition~5.2 in Cutkosky-Kurano~\cite{CK}).
 
\vspace{2mm}

In the rest of this section,
we state the results of this paper precisely.

Let $a$, $b$, $c$ be pairwise coprime integers.
We regard the polynomial ring $S = K[x, y, z]$ as a ${\Bbb Z}$-graded ring
by $\deg(x) = a$, $\deg(y) = b$ and $\deg(z) = c$.
Let ${\frak p}_K(a, b, c)$ be the kernel of the $K$-algebra
homomorphism
\[
\phi_K : S \longrightarrow K[t]
\]
given by $\phi_K(x) = t^a$, $\phi_K(y) = t^b$, $\phi_K(z) = t^c$.
If no confusion is possible, we simply denote ${\frak p}_K(a, b, c)$ by ${\frak p}$.

By a result of Herzog~\cite{Her}, we know that ${\frak p}_K(a, b, c)$
is generated by at most three binomials.
We define $s$, $t$, $u$ to be
\begin{equation}\label{Kstu}
\begin{array}{lll}
s = \min\{  {\Bbb N}a \cap ({\Bbb N}_0b+{\Bbb N}_0c) \}, \\
t = \min\{  {\Bbb N}b \cap ({\Bbb N}_0a+{\Bbb N}_0c) \}, \\
u = \min\{  {\Bbb N}c \cap ({\Bbb N}_0a+{\Bbb N}_0b) \},
\end{array}
\end{equation}
where ${\Bbb N}$ (resp.\ ${\Bbb N}_0$) denotes the set of positive integers
(resp.\ non-negative integers).
Let $t_1$, $u_1$, $s_2$, $u_2$, $s_3$, $t_3$ be non-negative integers such that
$sa = t_1b + u_1c$, $tb = s_2a + u_2c$, $uc = s_3a + t_3b$.
Then, ${\frak p}_K(a, b, c)$ is minimally generated by three elements
if and only if $s$, $t$, $u \ge 2$.
When this is the case, ${\frak p}_K(a, b, c)$ is minimally generated by three elements
\begin{equation}\label{Kseisei}
x^s - y^{t_1}z^{u_1}, \ \ 
y^t - x^{s_2}z^{u_2}, \ \ 
z^u - x^{s_3}y^{t_3},
\end{equation}
and $t_1$, $u_1$, $s_2$, $u_2$, $s_3$, $t_3$ must be positive integers
satisfying $s = s_2+s_3$, $t=t_1+t_3$, $u=u_1+u_2$.

For a prime ideal $P$ of $S$,
we define the symbolic Rees ring of $P$ to be
\[
\srees{P} = \oplus_{n \ge 0} P^{(n)}T^n \subset S[T] ,
\]
where $P^{(n)} = P^nS_P \cap S$ is the $n$th symbolic power of $P$
and $T$ is an indeterminate.
Here, $\srees{P}$ is a Noetherian ring if and only if 
$\srees{P}$ is finitely generated over $S$ as a ring.

In Section~\ref{Huneke's condition and mod $p$ reduction},
we give a simple proof to Huneke's criterion~\cite{Hu}.
We slightly generalize Huneke's criterion here.
Furthermore, we develop the method of mod $p$ reduction
introduced in Goto-Nishida-Watanabe~\cite{GNS}.

In Section~\ref{Ksec2}, we give a proof to the following theorem:

\begin{Theorem}\label{Kmain1}
Let $a$, $b$, $c$ be pairwise coprime positive integers.
Assume the following three conditions:
\begin{itemize}
\item[(i)] 
$K$ is a field of characteristic $0$,
\item[(ii)]
${\frak p}_K(a, b, c)$ is minimally generated by
the three elements as in (\ref{Kseisei}),
\item[(iii)]
$uc < \sqrt{abc}$.
\end{itemize}

Then, $\srees{{\frak p}_K(a, b, c)}$ is a Noetherian ring
if and only if there exists $\eta$ in $[{\frak p}^{(u)}]_{ab}$
such that $z^u -  x^{s_3}y^{t_3}$ and $\eta$ satisfy
Huneke's condition~\cite{Hu} (see Theorem~\ref{N2f}), that is,
\begin{equation}\label{Khuneke}
\length{S}{S/(x, z^u -  x^{s_3}y^{t_3}, \eta)} = u \cdot \length{S}{S/(x) + {\frak p}}
\end{equation}
holds.
\end{Theorem}

The condition (iii) as above implies that $z^u - x^{s_3}y^{u_3}$ is the negative curve,
that is, there exists the negative curve $C$ such that $C.E = 1$.
Theorem~\ref{Kmain1} says that there exists a curve $D$ such that $D \cap C = \emptyset$ and
$D \sim abA - uE$ if and only if $\srees{{\frak p}_\qq(a, b, c)}$ is a Noetherian ring,
where $A$ is an Weil divisor on $X$ satisfying ${\mathcal O}_X(A) = \pi^*{\mathcal O}_{\Bbb P}(1)$.

We emphasis that it is possible to verify weather there exists
$\eta$ in $[{\frak p}^{(u)}]_{ab}$ satisfying (\ref{Khuneke}) as above
using computers.
We shall prove this theorem using
the mod $p$ reduction as in Goto-Nishida-Watanabe~\cite{GNW},
Fujita's vanishing theorem (Theorem~1.4.35 in \cite{La}) and Cutkosky's method \cite{Cu} in characteristic $p > 0$.
The most important point is that the negative curve is isomorphic to 
${\Bbb P}^1_K$ in this case.

In Section~\ref{Ksec3}, we introduce the condition EU.
In Ebina~\cite{Ebina} and Uchisawa~\cite{Uchi},
the condition EU was defined and they proved that the condition EU is a sufficient condition
for finite generation under the assumptions (i), (ii), (iii) in Theorem~\ref{Kmain1}.
For the convenience of the reader, we shall give a proof of it in this paper.
Furthermore, in the case where $u \le 6$,
we show that the condition EU is a necessary and sufficient condition for the finite generation of the symbolic Rees ring of ${\frak p}$ in Propositin~\ref{Kmain2}.

In Section~\ref{Nsec4} we give an example of infinitely generated symbolic Rees ring of ${\frak p}$
where the homogeneous element of the minimal degree 
in  ${\frak p}^{(2)}$ is the negative curve in Example~\ref{N5f}.
We emphasis that 
one of the minimal generators of ${\frak p}$ is the negative curve in 
all known examples of infinitely generated $\srees{{\frak p}_K(a, b, c)}$,
except for this example.

\section{Huneke's condition and mod $p$ reduction}\label{Huneke's condition and mod $p$ reduction}

Let $S = K[x, y, z]$,
where $K$ is a field and $x, y, z$ are indeterminates.
We regard $S$ as a $\zz$-graded ring putting suitable weights on $x, y$ and $z$.
We set $\gm = (x, y, z)S$ and $R = S_\gm$.
Let $I$ be a homogeneous proper ideal of $S$
satisfying the following conditions;
\begin{itemize}
\item
$(x) + I$ is $\gm$-primary,

\item
$\ass_S S / I = \assh_S S / I :=
\{ \gp \in \ass_S S / I \mid 
\dim S / \gp = \dim S / I \}$, and

\item
$I_\gp$ is generated by $2$ elements for any $\gp \in \assh_S S / I$.
\end{itemize}
Then, $S / I$ is a $\zz$-graded Cohen-Macaulay ring of
$\dim S / I = 1$.
If we replace $x$ in the first assumption stated above with $y$ or $z$,
it can play the same role as $x$ in the arguments of this section.
So, homogeneous prime ideals of height $2$ are typical examples of $I$.
For any $n \in \zz$, we set
\[
I^{(n)} = \bigcap_{\gp \in \assh_S S / I}\,(I_\gp^n \cap S) ,
\]
where $I_\gp^n$ denotes the ideal $(I^n)_\gp = (I_\gp)^n$ of $S_\gp$.
Then, we have $\ass_S S / I^{(n)} = \assh_S S / I$ if $n > 0$,
and the equality
$I^{(n)} = I^n :_S x^i$ holds for $i \gg 0$,
which means that $I^{(n)}$ is a homogeneous ideal of $S$
and $(I^{(n)})_x = I_x^n$, where $I_x^n$ denotes the ideal 
$(I^n)_x = (I_x)^n$ of $S_x$.
Moreover, we set
\begin{eqnarray*}
\srees{I} & = & \sum_{n \geq 0}\, I^{(n)}T^n \subset S[ T ], \\
\serees{I} & = & \sum_{n \in \zz}\,I^{(n)} T^n \subset S[ T, T^{-1} ],\,\,\mbox{and} \\
\sgr{I} & = & \serees{I} / T^{-1}\serees{I} \hspace{1ex} = \hspace{1ex} 
\oplus_{n \geq 0}\, I^{(n)} / I^{(n + 1)},
\end{eqnarray*}
where $T$ is an indeterminate
and $I^{(n)} = S$ for $n \leq 0$.
Let us call $\srees{I}$ the symbolic Rees ring of $I$.

We set $\ga = I_\gm = IR$.
It is easy to see that $R / \ga$ is a Cohen-Macaulay local ring
of $\dim R / \ga = 1$ and
$\ass_R R / \ga = \assh_R R / \ga = \{ \gp R \mid \gp \in \assh_S S / I \}$.
Moreover, for $\gp \in \assh_S S / I$, we have $\ga_{\gp R} = I_\gp$,
which becomes a parameter ideal of $R_{\gp R} = S_\gp$.
For any $n \in \zz$, we set
\[
\ga^{(n)} = \bigcap_{P \in \assh_R R / \ga}\,(\ga_P^n \cap R) .
\]
Then, we have $\ga^{(n)} = I^{(n)}R$ and
$\ass_R R / \ga^{(n)} = \assh_R R / \ga$ if $n > 0$.
As $\ga^{(n)} = \ga^n :_R x^i$ holds for $i \gg 0$,
we have $(\ga^{(n)})_x = \ga_x^n$.
The $R$-algebras $\srees{\ga}$ and $\sgr{\ga}$ are derived from
$\srees{I}$ and $\sgr{I}$ respectively applying $R \otimes_S \ast$\,.
If $\srees{\ga}$ is finitely generated,
then there exists $0 < m \in \zz$ such that
$\ga^{(mn)} = (\ga^{(m)})^n$ for any $n \in \zz$.
This equality implies $I^{(mn)} = (I^{(m)})^n$ since
$I^{(mn)} \supseteq (I^{(m)})^n$ and $(I^{(m)})^n$ is a homogeneous ideal.
Thus we see that
$\srees{I}$ is finitely generated if so is $\srees{\ga}$.
The converse of this assertion holds obviously.

For a proper ideal $J$ of $S$ such that $S / J$ is Artinian,
we have $\length{S}{S / J} \geq \length{R}{R / JR}$,
and the equality holds if and only if $J$ is $\gm$-primary,
which holds if $J$ is homogeneous.

The purpose of this section is to review the condition on $I$
for its symbolic Rees ring to be finitely generated,
which was originally given by Huneke \cite{Hu}
in the case where $I$ is a prime ideal of a $3$-dimensional regular local ring.
Furthermore, using mod $p$ reduction technique for
prime numbers $p \gg 0$,
we give a condition on $I$ for
$\srees{I}$ to be infinitely generated,
which is a modification to the method introduced in \cite{GNW}.

Let us begin with the following

\begin{Proposition}\label{N2a}
Let $0 < k, \ell \in \zz$, $\xi \in I^{(k)}$ and $\eta \in I^{(\ell)}$.
Then we have
\[
\length{R}{R / (x, \xi, \eta)R} \geq
kl\cdot\length{S}{S / (x) + I},
\]
where the equality holds if and only if
$\ga \subseteq \sqrt{(\xi, \eta)R}$ and
\[
\length{S_\gp}{S_\gp / (\xi, \eta)S_\gp} = k\ell\cdot\length{S_\gp}{S_\gp / I_\gp}
\]
for all $\gp \in \assh_S S / I$.
\end{Proposition}

In order to prove Proposition~\ref{N2a}, let us recall the following fact.

\begin{Lemma}\label{N2c}
Let $A$ be a $2$-dimensional Cohen-Macaulay local ring and
$Q$ a parameter ideal of $A$.
Let $0 < k, \ell \in \zz$,
$\xi \in Q^k$ and $\eta \in Q^\ell$.
We assume that $\xi, \eta$ is an sop for $A$.
Then, we have
\[
\length{A}{A / (\xi, \eta)} \geq k\ell\cdot\length{A}{A / Q},
\]
and the equality holds if and only if one of the following conditions,
which are equivalent to each other,
is satisfied;
\begin{itemize}
\item[\rm (1)]
$QT \subseteq \sqrt{(\xi T^k, \eta T^\ell)\rees{Q}}$,
where $\rees{Q} = \sum_{n \geq 0} Q^nT^n \subset A[ T ]$,
\item[\rm (2)]
$\xi T^k, \eta T^\ell$ is an sop for $\gr{Q} = \rees{Q} / Q\rees{Q}$,
\item[\rm (3)]
$Q^{k + \ell - 1} = \xi Q^{\ell - 1} + \eta Q^{k - 1}$,
\item[\rm (4)]
$Q^n \cap (\xi, \eta)A = \xi Q^{n - k} + \eta Q^{n - \ell}$
for any $n \in \zz$.
\end{itemize}
\end{Lemma}

\noindent
{\it Proof.}
We set $J = (\xi^\ell, \eta^k)A$.
Then we have
\[
k\ell\cdot\length{A}{A / (\xi, \eta)} = \length{A}{A / J} = \mult{J}{A},
\]
where $\mult{J}{A}$ denotes the multiplicity of $A$ with respect to $J$.
Because $J \subseteq Q^{k\ell}$, it follows that
\[
\mult{J}{A} \geq \mult{Q^{k\ell}}{A} = (k\ell)^2\cdot\mult{Q}{A}
= (k\ell)^2\cdot\length{A}{A / Q}.
\]
Hence we get the required inequality.
Moreover, we see that the equality holds if and only if
$J$ is a reduction of $Q^{k\ell}$,
which is a condition equivalent to (1).
The equivalence of the conditions (1) and (2) is obvious.
Let us notice that $\gr{Q}$ is isomorphic to a polynomial ring
with $2$ variables over $A / Q$,
so its homogeneous sop is always a regular sequence,
which implies the equivalence of the conditions (2) and (4).
Moreover, if the condition (2) is satisfied,
it follows that $\gr{Q} / (\xi T^k, \eta T^\ell)\gr{Q}$ is an
Artinian $\zz$-graded ring whose a-invariant is $k + \ell - 2$ (cf. \cite{GW}),
so the equality of the condition (3) holds.
Finally, if the condition (3) is satisfied, we have
\[
(QT)^{k + \ell - 1} \subseteq 
\xi T^k\cdot Q^{\ell - 1}T^{\ell - 1} + \eta T^\ell\cdot Q^{k - 1}T^{k -1}
\subseteq (\xi T^k, \eta T^\ell)\rees{Q},
\]
and hence the condition (1) is satisfied.
\qed

\vspace{1em}
\noindent
{\it Proof of Proposition~\ref{N2a}.}
We may assume that $(x, \xi, \eta)R$ is $\gm R$-primary.
Then, as $R / (\xi, \eta)R$ is a Cohen-Macaulay $R$-module,
for which $x$ is an sop, we have
\[
\length{R}{R / (x, \xi, \eta)R} = \mult{xR}{R / (\xi, \eta)R}.
\]
Here we notice that
$\gp R \in \assh_R R / (\xi, \eta)R$
for any $\gp \in \assh_S S / I$.
Hence, using additive formula of multiplicity and Lemma~\ref{N2c}, we get
\begin{eqnarray*}
\mult{xR}{R / (\xi, \eta)R)} & = & \sum_{P \in \assh_R R / (\xi, \eta)R}
      \length{R_P}{R_P / (\xi, \eta)R_P}\cdot\mult{xR}{R / P}  \\
 & \geq & \sum_{\gp \in \assh_S S / I}
      \length{S_\gp}{S_\gp / (\xi, \eta)S_\gp}\cdot\mult{xR}{R / \gp R} \\
 & \geq & \sum_{\gp \in \assh_S S / I}
      k\ell\cdot\length{S_\gp}{S_\gp / I_\gp}\cdot\mult{xR}{R / \gp R} \\
 & = & k\ell\cdot\mult{xR}{R / \ga} \\
 & = & k\ell\cdot\length{R}{R / (x) + \ga} \\
 & = & k\ell\cdot\length{S}{S / (x) + I}.
\end{eqnarray*}
Thus we get the required inequality.
Moreover, we see that the equality holds if and only if
$\assh_R R / (\xi, \eta)R = \assh_R R / \ga$ and
$\length{S_\gp}{S_\gp / (\xi, \eta)S_\gp} = k\ell\cdot\length{S_\gp}{S_\gp / I_\gp}$
for all
$\gp \in \assh_S S / I$.
Since $\ga \subseteq \sqrt{(\xi, \eta)R}$ holds if and only if
$\assh_R R / (\xi, \eta)R = \assh_R R / \ga$,
the proof is complete.  
\qed   

\begin{Definition}\label{N2d}
Let $0 < k, \ell \in \zz, \xi \in I^{(k)}$ and $\eta \in I^{(\ell)}$.
We say that $\xi$ and $\eta$
satisfy Huneke's condition on $I$ (with respect to $x$) if
\[
\length{R}{R / (x, \xi, \eta)R} = k\ell\cdot\length{S}{S / (x) + I}.
\]
When this is the case,
for any $0 < i, j \in \zz$,
$\xi^i \in I^{(ki)}$ and $\eta^j \in I^{(\ell j)}$ also satisfy
Huneke's condition on $I$.
\end{Definition}

Even if there exist elements satisfying Huneke's condition,
those elements may not be homogeneous.
Although the existance of homogeneous elements
satisfying Huneke's condition is not clear,
but it can be verified elementary in special cases.
For example, the following remark implies that
if $\xi$ and $\eta$ satisfy Huneke's condition and
$\xi \equiv y^i$ mod $xS$ for some $0 < i \in \zz$,
then we can choose homogeneous parts of $\xi$ and $\eta$
so that they also satisfy Huneke's condition.

\begin{Lemma}\label{N2e}
Suppose $\xi \in S$ and $\xi \equiv y^i$ mod $xS$, where $0 < i \in \zz$.
Let $\xi'$ be the homogeneous part of $\xi$ containing $y^i$ as a term.
Then, the following assertions hold.
\begin{itemize}
\item[(1)]
$(x, \xi)S = (x, \xi')S = (x, y^i)S$.
\item[(2)]
For any $\eta \in S$, we can choose its homogeneous part $\eta'$ so that
\[
\length{S}{S / (x, \xi', \eta')} = \length{R}{R / (x, \xi, \eta)R}.
\]
\end{itemize}
\end{Lemma}

\noindent
{\it Proof.}
The assertion (1) holds obviously.
Let us verify the assertion (2).
We may assume $\eta \not\in (x, y)S$.
Then, as $x, y, \eta$ is an $R$-regular sequence, we have
\[
\length{R}{R / (x, \xi, \eta)R} = \length{R}{R / (x, y^i, \eta)R} =
i\cdot\length{R}{R / (x, y, \eta)R}.
\]
We write
\[
\eta \equiv \alpha_j z^j + \alpha_{j + 1} z^{j + 1} + \cdots \hspace{1ex}
\mbox{mod}\hspace{1ex} (x, y)S,
\]
where $0 < j \in \zz$ and $\alpha_j, \alpha_{j + 1}, \dots$
are elements of $K$ with $\alpha_j \neq 0$.
Since $\alpha_j + \alpha_{j + 1}z + \cdots$ is a unit of $R$, we have
\[
\length{R}{R / (x, y, \eta)R} = \length{R}{R / (x, y, z^j)R} = j.
\]
Thus we get
\[
\length{R}{R / (x, \xi, \eta)R} = ij.
\]
Let $\eta'$ be the homogeneous part of
$\eta$ containing $\alpha_j z^j$ as a term.
Then, as $\eta' \equiv \alpha_jz^j$ mod $(x, y)S$, it follows that
\[
\length{S}{S / (x, \xi', \eta')} = 
\length{S}{S / (x, y^i, \eta')} =
i\cdot\length{S}{S / (x, y, z^j)} = ij.
\]
Thus we get the required equality.
\qed

\begin{Theorem}\label{N2f}
The symbolic Rees algebra $\srees{I}$  is finitely generated over $R$
if and only if
there exist elements in $I^{(k)}$ and $I^{(\ell)}$
satisfying Huneke's condition on $I$ for some $0 < k, \ell \in \zz$.
\end{Theorem}

\noindent
{\it Proof.}
First, let us assume that
$\srees{I}$ is finitely generated.
Then, there exists a positive integer $m$ such that
$I^{(mn)} = (I^{(m)})^n$ for any $n \in \zz$.
We set $\gb = \ga^{(m)}$.
Then, for any $0 < n \in \zz$,
we have $\ga^{(mn)} = \gb^n$,
which means $\depth_R R / \gb^n = 1$.
Hence, by Burch's theorem (cf. \cite{Bu}), we see
\[
2 = \height_R \gb \leq \lambda(\gb) \leq 3 - 
\inf_{n > 0}\{ \depth R / \gb^n \} = 2,
\]
where $\lambda(\gb)$ denotes the Krull dimension of $R / \gm \otimes \gr{\gb}$,
which is called the analytic spread of $\gb$.
Thus we get $\lambda(\gb) = 2$.
Hence, we can choose $0 < i, j \in \zz$, $\xi \in I^{(mi)}$ and
$\eta \in I^{(mj)}$ such that $\xi T^i, \eta T^j$ is an sop for $R / \gm \otimes \gr{\gb}$.
(Here, we notice that we don't have to assume that $K$ is infinite
since we don't require $i = j = 1$.)
Let us take $r \gg 0$.
Then, we have $\gb^r = \xi\gb^{r - i} + \eta\gb^{r - j}$,
which means $\ga^{mr} \subseteq (\xi, \eta)R$, 
and so $\ga \subseteq \sqrt{(\xi, \eta)R}$.
Moreover, if $\gp \in \assh_S S / I$ and $mr < n$, we have
\[
I_\gp^n = \gb_{\gp R}^r I_\gp^{n - mr} =
(\xi\gb_{\gp R}^{r - i} + \eta\gb_{\gp R}^{r - j})I_\gp^{n - rm} =
\xi I_\gp^{n - mi} + \eta I_\gp^{n - mj},
\]
which means that $\xi T^{mi}, \eta T^{mj}$ is an sop for $\gr{I_\gp}$.
Therefore, by Proposition~\ref{N2a} and Lemma~\ref{N2c}, it follows that
$\xi$ and $\eta$ satisfy Huneke's condition on $I$.

Next, we assume that
there exist $0 < k, \ell \in \zz$, $\xi \in I^{(k)}$ and $\eta \in I^{(\ell)}$
such that $\xi$ and $\eta$
satisfy Huneke's condition on $I$.
We set $m = k\ell$, $\gb = \ga^{(m)}$ and 
$\gc = (\xi^\ell, \eta^k)R \subseteq \gb$.
Let us look at the exact sequence
\[
0 \ra \gc^r / \gb\gc^r \ra R / \gb\gc^r \ra R / \gc^r \ra 0
\]
of $R$-modules, where $r$ is any non-negative integer.
Since $\xi^\ell, \eta^k$ is an $R$-regular sequence,
$\gc^r / \gc^{r + 1}$ is $R / \gc$-free,
so $\gc^r / \gb\gc^r \cong R / \gb \otimes_R \gc^r / \gc^{r + 1}$
is $R / \gb$-free, which means
\[
\ass_R \gc^r / \gb\gc^r = \ass_R R / \gb = \assh_R R / \ga.
\]
On the other hand,
by Proposition~\ref{N2a} we have
\[
\ass_R R / \gc^r = \assh_R R / \gc = \assh_R R / \ga
\]
since $\xi^\ell$ and $\eta^k$ also satisfy Huneke's condition on $I$.
Thus we see
\[
\ass_R R / \gb\gc^r = \assh_R R / \ga.
\]
Now we take any $P \in \assh_R R / \ga$,
and write $P = \gp R$,
where $\gp \in \assh_S S / I$.
Then, by Proposition~\ref{N2a} and Lemma~\ref{N2c} we have
$I_\gp^{2m - 1} = \xi^\ell I_\gp^{m - 1} + \eta^k I_\gp^{m - 1}$,
which means $\gb_P^2 = (\gb\gc)_P$, and so $\gb_P^{r + 1} = (\gb\gc^r)_P$.
Hence, we get
\[
\ga^{(mr + m)} = \bigcap_{P \in \assh_R R / \ga} (\gb_P^{r + 1} \cap R) =
\bigcap_{P \in \assh_R R / \ga}((\gb\gc^r)_P \cap R) =
\gb\gc^r \subseteq \gb^{r + 1} \subseteq \ga^{(mr + m)},
\]
and hence $\ga^{(mr + m)} = \gb\gc^r = \gb^{r + 1}$.
Thus we see that the $m$-th veronese subring of
$\srees{\ga}$ is generated in degree one.
Therefore $\srees{\ga}$ is
Noetherian by \cite[Lemma (2.4)]{GHNV}.
Then $\srees{I}$ itself must be Noetherian.
\qed

\begin{Lemma}\label{N2g}
Let $0 < k, \ell \in \zz$, $\xi \in I^{(k)}$ and $\eta \in I^{(\ell)}$.
Suppose that $\xi$ and $\eta$ satisfy Huneke's condition on $I$.
Then the following assertions hold.
\begin{itemize}
\item[{\rm (1)}]
$\srees{\ga}_+ = \sqrt{(\xi T^k, \eta T^\ell)\srees{\ga}}$,
and hence
$\sgr{\ga}_+ \subseteq \sqrt{(\xi T^k, \eta T^\ell)\sgr{\ga}}$.
\item[{\rm (2)}]
$\ga^{(k + \ell - 1)} \subseteq (\xi, \eta)R$.
\item[{\rm (3)}]
$\ga_x^n \cap (\xi, \eta)R_x = \xi\ga_x^{n - k} + \eta\ga_x^{n - \ell}$
for any $n \in \zz$.
\item[{\rm (4)}]
$\ga^{(n)} \cap (\xi, \eta)R = \xi\ga^{(n - k)} + \eta\ga^{(n - \ell)}$ if $n \leq k + \ell$.
\item[{\rm (5)}]
If $k = 1$ or $2$, then we have
\[
\ga^{(n)} \cap (\xi, \eta)R = \xi\ga^{(n - k)} + \eta\ga^{(n - \ell)}
\]
for any $n \in \zz$, which means that
$\xi T^k, \eta T^\ell$ is a regular sequence on $\sgr{\ga}$,
and hence $\grade \sgr{\ga}_+ = 2$.
\end{itemize}
\end{Lemma}

\noindent
{\it Proof.}
(1) We set $m = k\ell, \gb = \ga^{(m)}$ and $\gc = (\xi^\ell, \eta^k)R$.
Then, as is stated in the proof of Theorem~\ref{N2f},
we have $\ga^{(mr + m)} = \gb^{r +1} = \gb\gc^r$
for any $0 \leq r \in \zz$.
Let us take any $0 < n \in \zz$ and $\rho \in \ga^{(n)}$.
Then we have $\rho^{2m} \in \ga^{(m(2n - 1) + m)} = \gb^{2n}
= \gb\gc^{2n - 1} \subseteq \gc\gb^{2n - 1} =
\xi^\ell \gb^{2n - 1} + \eta^k \gb^{2n - 1} \subseteq
\xi \ga^{(2mn - k)} + \eta \ga^{(2mn - \ell)}$, so
\[
(\rho T^n)^{2m} \in \xi T^k\cdot \ga^{(2mn - k)}T^{2mn - k} + 
       \eta T^\ell\cdot \ga^{(2mn - \ell)}T^{2mn - \ell}.
\]
Hence we get the assertion (1)

(2) 
Let us take any $P \in \assh_R R / \ga$ and write $P = \gp R$,
where $\gp \in \assh_S S / I$.
Then, as $R_P = S_\gp$ and $\ga_P = I_\gp$,
by Proposition~\ref{N2a} and Lemma~\ref{N2c}, we have
$\ga_P^{k + \ell - 1} = \xi\ga_{P}^{\ell - 1} + \eta\ga_P^{k - 1}
\subseteq (\xi, \eta)R_P$.
Therefore we get
\[
\ga^{(k + \ell - 1)} = 
\bigcap_{P \in \assh_R R / \ga}(\ga_P^{k + \ell - 1} \cap R) \subseteq
\bigcap_{P \in \assh_R R / \ga}((\xi, \eta)R_P \cap R) =
(\xi, \eta)R.
\]

(3)
Since $\xi \in I_x^k$ and $\eta \in I_x^\ell$,
the inclusion $\ga_x^n \supseteq \xi\ga_x^{n - k} + \eta\ga_x^{n - \ell}$ holds obviously.
So, it is enough to show
\[
\ga_P^n \cap (\xi, \eta)R_P = \xi\ga_P^{n - k} + \eta\ga_P^{n - \ell}
\]
for any $P \in \spec R$ satisfying
$\xi\ga^{n - k} + \eta\ga^{n - \ell} \subseteq P$ and $x \not\in P$.
Such $P$ must contains $\ga$ since $\ga \subseteq \sqrt{(\xi, \eta)R}$,
so there exists $\gp \in \assh_S S / I$ such that $P = \gp R$.
Then, by Proposition~\ref{N2a} and Lemma~\ref{N2c}, we get the required equality
as $R_P = S_\gp$ and $\ga_P = I_\gp$.

(4) 
Let $n \leq k + \ell$ and $\varphi \in \ga^{(n)} \cap (\xi, \eta)R$.
We write $\varphi = \xi u + \eta v$, where $u, v \in R$.
Since $\varphi \in \ga_x^n \cap (\xi, \eta)R_x = 
\xi\ga_x^{n - k} + \eta\ga_x^{n - \ell}$ by (3),
there exist $\alpha \in \ga_x^{n - k}$ and $\beta \in \ga_x^{n - \ell}$
such that $\varphi = \xi\alpha + \eta\beta$.
Here, we take $i \gg 0$ so that $x^i\alpha \in \ga^{n - k}$ and $x^i\beta \in \ga^{n - \ell}$.
Then we have
$x^i(\xi u + \eta v) = x^i\varphi = x^i(\xi\alpha + \eta\beta)$, so
$\xi(x^i u - x^i\alpha) = \eta(x^i\beta - x^i v)$.
Since $\xi, \eta$ is an $R$-regular sequence,
it follows that
$x^i u - x^i\alpha \in \eta R \subseteq \ga^{(\ell)} \subseteq \ga^{(n - k)}$ and
$x^i\beta - x^i v \in \xi R \subseteq \ga^{(k)} \subseteq \ga^{(n - \ell)}$.
Hence $x^i u \in \ga^{(n - k)}$ and $x^i v \in \ga^{(n - \ell)}$,
which means $u \in \ga^{(n - k)}$ and $v \in \ga^{(n - \ell)}$.
Thus we get $\varphi \in \xi\ga^{(n - k)} + \eta\ga^{(n - \ell)}$.

(5)
Let $k = 1$ or $2$.
By (2) and (4), it is enough to show 
\[
\ga^{(n)} = \xi\ga^{(n - k)} + \eta\ga^{(n - \ell)}
\]
assuming $n > k + \ell$.
We take positive integers $m$ and $r$ such that
$n - \ell = km - r$ and $0 \leq r < k$.
Then, $m \geq 2$ and $r$ is $0$ or $1$.
Since $\xi^m \in I^{(km)}$ and $\eta \in I^{(\ell)}$
also satisfy Huneke's condition on $I$ and
$km + \ell - 1 \leq km + \ell - r = n \leq km + \ell$,
we have $\ga^{(n)} \subseteq (\xi^m, \eta)R$ by (2) and
$\ga^{(n)} \cap (\xi^m, \eta)R = \xi^m\ga^{(\ell - r)} + \eta\ga^{(n - \ell)}$ by (4).
Let us notice that $\xi^{m - 1}\ga^{(\ell - r)} \subseteq \ga^{(n - k)}$ as
$k(m - 1) + (\ell - r) = n - k$.
Thus we see $\ga^{(n)} \subseteq \xi\ga^{(n - k)} + \eta\ga^{(n - \ell)}$.
Since the converse inclusion is obvious, we get the required equality.
\qed

\begin{Definition}\label{N2h}
Let $0 < k \in \zz$ and $\xi \in I^{(k)}$.
We denote by $\hc{I}{k}{\xi}$ the set of positive integers $\ell$
for which there exists $\eta \in I^{(\ell)}$ such that
$\xi$ and $\eta$ satisfy Huneke's condition on $I$.
\end{Definition}

\begin{Remark}\label{N2o}
Let $k$ and $\xi$ be as in Definition \ref{N2h}.
If $\xi \equiv y^i$ mod $xS$, where $0 < i \in \zz$,
and $\xi'$ is the homogeneous part of $\xi$ containing $y^i$ as a term,
we have $\hc{I}{k}{\xi} = \hc{I}{k}{\xi'}$ by Lemma \ref{N2e} (1).
\end{Remark}

\begin{Proposition}\label{N2i}
Let $k = 1$ or $2$, and let $\xi \in I^{(k)}$.
Suppose that $\xi \equiv y^i \,\,{\rm mod}\,\, xS$
for some $0 < i \in \zz$ and $\hc{I}{k}{\xi} \neq \phi$.
We set $m = \min\hc{I}{k}{\xi}$.
Then the following assertions hold.
\begin{itemize}
\item[{\rm (1)}]
$\hc{I}{k}{\xi} = \{ m, 2m, 3m, \cdots \}$.
\item[{\rm (2)}]
$S[\{ I^{(n)}T^n \mid 1 \leq n \leq m - 1 \}] \subsetneq \srees{I}$.
\item[{\rm (3)}]
If there exist elements in $I^{(k')}$ and $I^{(\ell')}$
satisfying Huneke's condition on $I$ for $0 < k', \ell' \in \zz$,
we have
\[
S[\{ I^{(n)}T^n \mid 1 \leq n \leq \max\{ k', \ell', k' + \ell' - 2 \} \}] = \srees{I}.
\]
In particular, 
\[
S[\{ I^{(n)}T^n \mid 1 \leq n \leq \max\{ k, m \} \}] = \srees{I}.
\]
\end{itemize}
\end{Proposition}

\noindent
{\it Proof.}
By Remark \ref{N2o}, 
we may assume that $\xi$ is homogeneous.
Then, by Lemma \ref{N2e} (2), we can choose a homogeneous element $\eta \in I^{(m)}$
such that $\xi$ and $\eta$ satisfy Huneke's condition on $I$.

(1)
We obviously have $\hc{I}{k}{\xi} \supseteq \{ m, 2m, 3m, \cdots \}$.
In order to show the converse inclusion,
We suppose that there exists $\ell \in \hc{I}{k}{\xi}$ which is not a multiple of $m$.
Let us choose such $\ell$ as small as possible.
Then, there exists a homogeneous element $\rho \in I^{(\ell)}$
such that $\xi$ and $\rho$ satisfy Huneke's condition on $I$.
Since $m < \ell$, by Lemma~\ref{N2g} (2) and (5), we have
$\ga^{(\ell)} = \xi\ga^{(\ell - k)} + \eta\ga^{(\ell - m)}$, which implies
\[
I^{(\ell)} = \xi I^{(\ell - k)} + \eta I^{(\ell - m)}
\]
as $\xi$ and $\eta$ are homogeneous.
Hence, there exists a homogeneous element $\rho' \in I^{(\ell - m)}$ such that
\[
\rho \equiv \eta\rho' \hspace{1ex}\mbox{mod}\hspace{1ex} \xi I^{(\ell - k)}.
\]
Then $\rho \in (\xi, \rho')S$, and hence we get
\[
\ga \subseteq \sqrt{(\xi, \rho')R}
\]
as $\ga \subseteq \sqrt{(\xi, \rho)R}$ by Proposition~\ref{N2a}.
Now we take any $\gp \in \assh_S S / I$ and $n \gg 0$.
Then, by Proposition~\ref{N2a} and Lemma~\ref{N2c}, we have
\[
I_\gp^n = \xi I_\gp^{n - k} + \rho I_\gp^{n - \ell} =
\xi I_\gp^{n - k} + \eta\rho' I_\gp^{n - \ell} \subseteq
\xi I_\gp^{n - k} + \rho' I_\gp^{n - (\ell - m)} \subseteq I_\gp^n,
\]
so we get
\[
I_\gp^n = \xi I_\gp^{n - k} + \rho' I_\gp^{n - (\ell-m)}.
\]
Therefore, $\xi$ and $\rho'$ satisfy Huneke's condition on $I$,
so $\ell - m \in \hc{I}{k}{\xi}$, which contradicts to the minimality
of $\ell$ as $\ell - m$ is not a multiple of $m$.
Consequently, we see that any $\ell \in \hc{I}{k}{\xi}$ is a multiple of $m$.

(2)
The assertion holds obviously if $m = 1$,
so let us consider the case where $m \geq 2$.
Suppose
\[ 
\eta T^m \in S[\{ I^{(n)}T^n \mid 1 \leq n \leq m - 1 \}].
\]
Then we have
\[
\eta \in \sum_{\alpha = 1}^{m - 1}\,I^{(\alpha)} I^{(m - \alpha)}.
\]
We set $\overline{S} = S / (x, y) \cong K[z]$.
Since any homogeneous ideal of $\overline{S}$ is a power of $z\overline{S}$,
\[
\sum_{\alpha = 1}^{m - 1}\,I^{(\alpha)}I^{(m - \alpha)}\overline{S} =
I^{(\beta)}I^{(m - \beta)}\overline{S}
\]
holds for some $\beta = 1, 2, \dots, m - 1$.
Moreover, we can choose homogeneous elements
$\rho \in I^{(\beta)}$ and $\rho' \in I^{(m - \beta)}$
such that $\eta$ and $\rho\rho'$ have the same class in $\overline{S}$,
which is equivalent to
\[
\eta \equiv \rho\rho' \hspace{1ex}\mbox{mod}\hspace{1ex} (x, y).
\]
Then, by Proposition~\ref{N2a}, we have
\begin{eqnarray*}
\length{S}{S / (x, \xi, \rho)} & \geq & 
    k\beta\cdot\length{S}{S / (x) + I} \hspace{2ex}\mbox{and} \\
\length{S}{S / (x, \xi, \rho')} & \geq & 
    k(m - \beta)\cdot\length{S}{S / (x) + I}.
\end{eqnarray*}
Since $(x, y, \eta)$, $(x, y, \rho)$ and $(x, y, \rho')$ are all
homogeneous $\gm$-primary ideals, we have
\begin{eqnarray*}
\length{S}{S / (x, \xi, \eta)} & = &
          \length{S}{S / (x, y^i, \eta)}  \\
 & = & i\cdot\length{S}{S / (x, y, \eta)}  \\
 & = & i\cdot\length{S}{S / (x, y, \rho\rho')}  \\
 & = & i\cdot\{\length{S}{S / (x, y, \rho)} + \length{S}{S / (x, y, \rho')}\} \\
 & = & \length{S}{S / (x, y^i, \rho)} + \length{S}{S / (x, y^i, \rho')}  \\
 & = & \length{S}{S / (x, \xi, \rho)} + \length{S}{S / (x, \xi, \rho')} \\
 & \geq & k\beta\cdot\length{S}{S / (x) + I} +
                  k(m - \beta)\cdot\length{S}{S / (x) + I}  \\
 & = & \{ k\beta + k(m - \beta) \}\cdot\length{S}{S / (x) + I} \\
 & = & km\cdot\length{S}{S / (x) + I}  \\
 & = & \length{S}{S / (x, \xi, \eta)}.
\end{eqnarray*}
Consequently, it follows that
\begin{eqnarray*}
\length{S}{S / (x, \xi, \rho)} & = & 
    k\beta\cdot\length{S}{S / (x) + I} \hspace{2ex}\mbox{and} \\
\length{S}{S / (x, \xi, \rho')} & = & 
    k(m - \beta)\cdot\length{S}{S / (x) + I}.
\end{eqnarray*}
Hence we get $\beta, m - \beta \in \hc{I}{k}{\xi}$,
which contradicts to the minimality of $m$.
Thus we see
\[
\eta T^m \not\in S[\{ I^{(n)}T^n \mid 1 \leq n \leq m - 1 \}].
\]

(3)
Let $0 < k', \ell' \in \zz$, $\xi' \in I^{(k')}$ and $\eta' \in I^{(\ell')}$.
Suppose that $\xi'$ and $\eta'$ satisfy Huneke's condition on $I$.
Then, by Lemma~\ref{N2g} (1), we have
\[
\sgr{\ga}_+ \subseteq \sqrt{(\xi' T^{k'}, \eta' T^{\ell'})\sgr{\ga}}.
\]
On the other hand, from the existance of $\xi$ and $\eta$, we see
$\grade \sgr{\ga}_+ = 2$ by Lemma~\ref{N2g} (5).
Hence, it follows that $\xi' T^{k'}, \eta' T^{\ell'}$
is a regular sequence on $\sgr{\ga}$.
If $k' + \ell' -1 \leq n$, we have $\ga^{(n)} \subseteq (\xi', \eta')R$
by Lemma~\ref{N2g} (2), so
\[
\ga^{(n)} = \ga^{(n)} \cap (\xi', \eta')R = \xi'\ga^{(n - k')} + \eta'\ga^{(n - \ell')}.
\]
Thus we see
\begin{eqnarray*}
\srees{\ga} & = & 
   S[ \xi' T^{k'}, \eta' T^{\ell'}, \{ \ga^{(n)}T^n \mid 1 \leq n \leq k' + \ell' -2 \} ] \\
 & = & S[ \{\ga^{(n)}T^n \mid 1 \leq n \leq \max\{ k', \ell', k' + \ell' - 2 \} \} ],
 \end{eqnarray*}
 which means that the first assertion of (3) holds.
 We get the last assertion taking $k$ and $m$ as $k'$ and $\ell'$, respectively.
\qed

\vspace{1em}
In the rest of this section,
let $S_{\zz} = \zz[x, y, z]$.
Moreover, for a field $K$,
we denote $K[x, y, z]$ by $S_K$ instead of $S$
in order to emphasize that the coefficient field is $K$.
Putting suitable weights on $x, y$ and $z$,
we regard $S_{\zz}$ and $S_K$ as $\zz$-graded rings.
We set $\gm_\zz = (x, y, z)S_{\zz}$, $\gm_K = (x, y, z)S_K$ and $R_K = (S_K)_{m_K}$.
When we denote an ideal of $S_\zz$ by $J_\zz$,
the ideal $J_\zz S_K$ is denoted by $J_K$.
Similarly, when we denote an element of $S_\zz$ by $\xi_\zz$,
its image in $S_K$ is denoted by $\xi_K$.
For a prime number $p$, we set $\fp = \zz / p\zz$.
Of course, $S_\qq = (\zz\setminus \{ 0 \})^{-1}S_\zz$ and
$S_\fp = S_\zz / pS_\zz$.

\begin{Lemma}\label{N2m}
Let $J_\zz$ be an ideal of $S_\zz$.
Then, we have
\[
\length{R_\qq}{R_\qq / (J_\qq)_{\gm_\qq}} =
\length{R_\fp}{R_\fp / (J_\fp)_{\gm_\fp}}
\]
for any prime number $p \gg 0$.
If $J_\zz$ is homogeneous, we may replace
$R_\qq,  (J_\qq)_{\gm_\qq}, R_\fp$ and $(J_\fp)_{\gm_\fp}$  with
$S_\qq, J_\qq, S_\fp$ and $J_\fp$, respectively.
\end{Lemma}

\noindent
{\it Proof.}
First, let us consider the case where
$R_\qq / (J_\qq)_{\gm_\qq}$ is Artinian.
We prove the required equality
by induction on $\length{R_\qq}{R_\qq / (J_\qq)_{\gm_\qq}}$.

If $\length{R_\qq}{R_\qq / (J_\qq)_{\gm_\qq}} = 0$, then
$J_\qq$ contains an element which does not belong to $\gm_\qq$,
so there exists $\xi_\zz \in J_\zz \setminus \gm_\zz$.
Let us take a prime number $p \gg 0$ so that
the constant term of $\xi_\zz$, which is non-zero,
is not a multiple of $p$.
Then, $\xi_{\fp} \in J_\fp \setminus \gm_{\fp}$.
Hence $(J_\fp)_{\gm_\fp} = R_\fp$, so 
$\length{R_\fp}{R_\fp / (J_\fp)_{\gm_\fp}} = 0$.

Now we suppose $\length{R_\qq}{R_\qq / (J_\qq)_{\gm_\qq}} > 0$.
Then, as $\gm_\zz \in \ass_{S_\zz} S_\zz / J_\zz$,
there exists $\eta_\zz \in S_\zz$ such that $J_\zz : \eta_\zz = \gm_\zz$.
We set $L_\zz = J_\zz + (\eta_\zz)$.
Let us notice $L_\qq / J_\qq \cong S_\qq / \gm_\qq$.
Hence we have $\length{R_\qq}{(L_\qq / J_\qq)_{\gm_\qq}} = 1$, so
\[
\length{R_\qq}{R_\qq / (L_\qq)_{\gm_\qq}} =
\length{R_\qq}{R_\qq / (J_\qq)_{\gm_\qq}} - 1.
\]
Here, we take a prime number $p \gg 0$.
Then the hypothesis of induction implies
\[
\length{R_\qq}{R_\qq / (L_\qq)_{\gm_\qq}} =
\length{R_\fp}{R_\fp / (L_\fp)_{\gm_\fp}}.
\]
Moreover, by taking larger $p$ if necessary,
we may assume that $p$ is regular on $S_\zz / L_\zz$.
If $\eta_\fp \in J_\fp$,
we have $\eta_\zz \in J_\zz + pS_\zz$,
so there exists $\rho_\zz \in S_\zz$ such that
$\eta_\zz \equiv p\cdot\rho_\zz$ mod $J_\zz$.
Since $p$ is regular on $S_\zz / L_\zz$, we have $\rho_\zz \in L_\zz$,
so there exists $\sigma_\zz \in S_\zz$ such that
$\rho_\zz \equiv \eta_\zz\sigma_\zz$ mod $J_\zz$.
Then, we have $\eta_\zz \equiv p\cdot\eta_\zz\sigma_\zz$ mod $J_\zz$,
and hence $1 - p\cdot\sigma_\zz \in J_\zz : \eta_\zz = \gm_\zz$,
which is impossible.
Thus we see $\eta_\fp \not\in J_\fp$.
Hence we have $\gm_\fp = J_\fp : \eta_\fp$ since
$\eta_\fp\cdot\gm_\fp \subseteq J_\fp$ holds obviously.
Then, we get $L_\fp / J_\fp \cong R_\fp / \gm_\fp$,
so $\length{R_\fp}{L_\fp / J_\fp} = 1$.
Consequently, 
\[
\length{R_\fp}{R_\fp / (L_\fp)_{\gm_\fp}} = 
\length{R_\fp}{R_\fp / (J_\fp)_{\gm_\fp}} - 1.
\]
Therefore, the required equality follows.

Next, we assume $\dim R_\qq / (J_\qq)_{\gm_\qq} > 0$,
and aim to prove $\dim R_\fp / (J_\fp)_{\gm_\fp} > 0$ for $p \gg 0$.
In this case, there exists $P_\zz \in \spec S_\zz$
such that $J_\zz \subseteq P_\zz \subsetneq \gm_\zz$.
Let us take any $\tau_\zz \in \gm_\zz \setminus P_\zz$
and choose a prime number $p \gg 0$ so that
$p$ is regular on $S_\zz / (\tau_\zz) + P_\zz$.
Then, as $p, \tau_\zz$ is a regular sequence on $(S_\zz / P_\zz)_{(pS_\zz + \gm_\zz)}$,
it follows that $\tau_\zz$ is regular on 
$(S_\zz / pS_\zz + P_\zz)_{(pS_\zz + \gm_\zz)} \cong R_\fp / (P_\fp)_{\gm_\fp}$.
Hence, we have $\dim R_\fp / (P_\fp)_{\gm_\fp} > 0$,
and so $\dim R_\fp / (J_\fp)_{\gm_\fp} > 0$ as $J_\fp \subseteq P_\fp$.
\qed

\vspace{1em}
In the rest of this section, let $I_\zz$ be a homogeneous ideal of $S_\zz$ contained in $\gm_\zz$.
We assume that the following conditions are satisfied for any field $K$;

\begin{itemize}
\item
$(x) + I_K$ is $\gm_K$-primary,
\item
$\ass_{S_K} S_K / I_K = \assh_{S_K} S_K / I_K$, and
\item
$(I_K)_\gp$ is generated by $2$ elements for any $\gp \in \assh_{S_K} S_K / I_K$.
\end{itemize}
Furthermore, for any $n \in \zz$, we set
$(I^{(n)})_\zz = \bigcup_{i > 0} ((I_\zz)^n :_{S_\zz} x^i)$,
which is a homogeneous ideal of $S_\zz$.
Let us denote $(I^{(n)})_\zz S_K$ by $(I^{(n)})_K$ for any field $K$.

\begin{Lemma}\label{N2n}
The following assertions hold for any $n \in \zz$.
\begin{itemize}
\item[{\rm (1)}]
$(I_\qq)^{(n)} = (I^{(n)})_\qq$.
\item[{\rm (2)}]
$(I_\fp)^{(n)} = (I^{(n)})_\fp$ for any prime number $p \gg 0$.
\end{itemize}
\end{Lemma}

\noindent
{\it Proof.}
First, let us notice that, for any field $K$,
we have $(I_K)^{(n)} = \bigcup_{i > 0} {((I_K)^n :_{S_K} x^i)}$,
and hence $(I_K)^{(n)} \supseteq (I^{(n)})_K$ holds.
The converse inclusion holds obviously if $K = \qq$.
So, we have to prove $(I_\fp)^{(n)} \subseteq (I^{(n)})_\fp$ for $p \gg 0$.

Let us take a prime number $p \gg 0$ so that
$p$ is regular on $S_\zz / (x) + (I^{(n)})_\zz$.
Moreover, we take any $\xi_\zz \in S_\zz$ satisfying 
$\xi_\fp \in (I_\fp)^{(n)}$.
Then, there exists $0 < i \in \zz$ such that
$x^i\xi_\fp \in (I_\fp)^n$,
which means $x^i\xi_\zz \in pS_\zz + (I_\zz)^n \subseteq pS_\zz + (I^{(n)})_\zz$.
Hence, there exists $\eta_\zz \in S_\zz$ such that
$x^i\xi_\zz \equiv p\eta_\zz$ mod $(I^{(n)})_\zz$.
Since $x, p$ is a regular sequence on $S_\zz / (I^{(n)})_\zz$,
so is $x^i, p$.
Hence $\eta_\zz \in (x^i) + (I^{(n)})_\zz$, so there exists
$\rho_\zz \in S_\zz$ such that
$\eta_\zz \equiv x^i\rho_\zz$ mod $(I^{(n)})_\zz$.
Then we have $x^i\xi_\zz \equiv px^i\rho_\zz$ mod $(I^{(n)})_\zz$,
which means $\xi_\zz - p\rho_\zz \in (I^{(n)})_\zz$.
Thus we get $\xi_\fp \in (I^{(n)})_\fp$.
\qed

\begin{Proposition}\label{N2k}
Let $0 < k, \ell \in \zz, \xi_\zz \in (I^{(k)})_\zz$ and
$\eta_\zz \in (I^{(\ell)})_\zz$.
Suppose that $\xi_\qq \in (I_\qq)^{(k)}$ and $\eta_\qq \in (I_\qq)^{(\ell)}$
satisfy Huneke's condition on $I_\qq$.
Then, for any prime number $p \gg 0$,
$\xi_\fp \in (I_\fp)^{(k)}$, $\eta_\fp \in (I_\fp)^{(\ell)}$,
and these elements satisfy Huneke's condition on $I_\fp$.
\end{Proposition}

\noindent
{\it Proof.}
Let $p \gg 0$.
Then, by Lemma \ref{N2n}, we have
$\xi_\fp \in (I_\fp)^{(k)}$ and $\eta_\fp \in (I_\fp)^{(\ell)}$.
Moreover, by Lemma \ref{N2m}, we have
\begin{eqnarray*}
\length{R_\fp}{R_\fp / (x, \xi_\fp, \eta_\fp)} & = &
  \length{R_\qq}{R_\qq / (x, \xi_\qq, \eta_\qq)}  \\
 & = & k\ell\cdot\length{S_\qq}{S_\qq / (x) + I_\qq} \\
 & = & k\ell\cdot\length{S_\fp}{S_\fp / (x) + I_\fp}.
\end{eqnarray*}
Thus we get the required assertion.
\qed

\begin{Theorem}\label{N2l}
Let $k = 1$ or $2$.
Let $\xi_\zz \in (I^{(k)})_\zz$ and
$\xi_\zz \equiv y^{i} \,\,{\rm mod}\,\, xS_\zz$
for some $0 < i \in \zz$.
Suppose that there exists a positive integer $r$
such that, for any prime number $p \gg 0$, we have
$rp^{e_p} \in \hc{I_\fp}{k}{\xi_\fp}$
for some $0 < e_p \in \zz$.
Then the following conditions are equivalent:
\begin{itemize}
\item[(1)]
$\srees{I_\qq}$ is finitely generated.
\item[(2)]
$\hc{I_\qq}{k}{\xi_\qq} \neq \emptyset$.
\item[(3)]
$r \in \hc{I_\qq}{k}{\xi_\qq}$.
\item[(4)]
$r \in \hc{I_\fp}{k}{\xi_\fp}$ for any prime number $p \gg 0$.
\end{itemize}
\end{Theorem}

\noindent
{\it Proof.}\,
Let $\xi'_\zz$ be the homogeneous part of $\xi_\zz$
containing $y^i$ as a term.
Then, as $\xi'_\zz \in (I^{(k)})_\zz$,
we have $\xi'_\qq \in (I_\qq)^{(k)}$ and $\xi'_\fp \in (I_\fp)^{(k)}$ for $p \gg 0$
by Lemma \ref{N2n}.
Moreover, by Remark \ref{N2o}, we have 
$\hc{I_\qq}{k}{\xi_\qq} = \hc{I_\qq}{k}{\xi'_\qq}$ and
$\hc{I_\fp}{k}{\xi_\fp} = \hc{I_\fp}{k}{\xi'_\fp}$ for $p \gg 0$.
Hence, by replacing $\xi_\zz$ with $\xi'_\zz$, 
we may assume that $\xi_\zz$ is homogeneous from the beginning.
It is easy to see $(3) \Rightarrow (2) \Rightarrow (1)$.

Now, we start to prove $(1) \Rightarrow (4)$.
By Theorem~\ref{N2f} and Lemma~\ref{N2n} (1),
there exist $0 < k', \ell' \in \zz, \zeta_\zz \in (I^{(k')})_\zz$
and $\rho_\zz \in (I^{(\ell')})_\zz$ such that
$\zeta_\qq \in (I_\qq)^{(k')}$ and $\rho_\qq \in (I_\qq)^{(\ell')}$
satisfy Huneke's condition on $I_\qq$.
Here, we take a prime number $p \gg 0$ such that
$\zeta_\fp \in (I_\fp)^{(k')}$, $\rho_\fp \in (I_\fp)^{(\ell')}$,
and these elements satisfy Huneke's condition on $I_\fp$,
which is possible by Proposition~\ref{N2k}.
By taking larger $p$ if necessary, we may assume $p > \max\{ k', \ell', k' + \ell' -2 \}$
and our assumption on $\hc{I_\fp}{k}{\xi_\fp}$ is satisfied.
Then, as $\hc{I_\fp}{k}{\xi_\fp} \neq \phi$, we have
\[
S_\fp[ \{ (I_\fp)^{(n)}T^n \mid 1 \leq n \leq p - 1 \} ] = \srees{I_\fp}
\]
by Proposition~\ref{N2i} (3).
We set $m = \min\hc{I_\fp}{k}{\xi_\fp}$ and take $0 < e_p \in \zz$
such that $rp^{e_p} \in \hc{I_\fp}{k}{\xi_\fp}$.
Then, by Proposition~\ref{N2i} (1), there exists $m' \in \zz$
such that $rp^{e_p} = mm'$.
Since Proposition~\ref{N2i} (2) implies $m < p$,
$m$ is not a multiple of $p$,
so $m'$ is a multiple of $p^{e_p}$.
Hence $r$ is a multiple of $m$, 
which means $r \in \hc{I_\fp}{k}{\xi_\fp}$.

Next, we shall prove $(4) \Rightarrow (3)$.
Let us take a prime number $p \gg 0$ such that
$r \in \hc{I_\fp}{k}{\xi_\fp}$,
$\length{S_\qq}{S_\qq / (x) + I_\qq} = \length{S_\fp}{S_\fp / (x) + I_\fp}$
and $(I_\fp)^{(r)} = (I^{(r)})_\fp$,
which is possible by Lemma~\ref{N2m} and Lemma~\ref{N2n}.
Then, by Lemma \ref{N2e} (2) and Lemma \ref{N2n}, 
there exists a homogeneous element  $\eta_\zz \in (I^{(r)})_\zz$
such that
$\xi_\fp \in (I_\fp)^{(k)}$ and $\eta_\fp \in (I_\fp)^{(r)}$
satisfy Huneke's condition on $I_\fp$.
We write
\[
\eta_\zz \equiv \alpha z^j
\hspace{2ex}
\mbox{mod}
\hspace{1ex}
(x, y)S_\zz,
\]
where $j$ is a positive integer and $\alpha$ is an integer
which is not a multiple of $p$.
Let $K = \qq$ or $\fp$.
Then, as the image of
$\alpha$ in $K$ is not vanished,
we have $(x, y, \eta_K)S_K = (x, y, z^j)S_K$.
Hence we get
\begin{eqnarray*}
\length{S_K}{S_K / (x, \xi_K, \eta_K)} & = & \length{S_K}{S_K / (x, y^i, \eta_K)} \\
  & = & i \cdot \length{S_K}{S_K / (x, y, \eta_K)} \\
  & = & i \cdot \length{S_K}{S_K / (x, y, z^j)} \\
  & = & ij.
\end{eqnarray*}
Consequently, we have
\begin{eqnarray*}
\length{S_\qq}{S_\qq / (x, \xi_\qq, \eta_\qq)} & = 
                 & \length{S_\fp}{S_\fp / (x, \xi_\fp, \eta_\fp)}  \\
  & = & kr\cdot\length{S_\fp}{S_\fp / (x) + I_\fp}  \\
  & = & kr\cdot\length{S_\qq}{S_\qq / (x) + I_\qq},
\end{eqnarray*}
which means $r \in \hc{I_\qq}{k}{\xi_\qq}$.
\qed

\section{Proof of Theorem~\ref{Kmain1}}\label{Ksec2}

We shall prove Theorem~\ref{Kmain1} in this section.

Let $K$ be a field and $a$, $b$, $c$ be pairwise coprime positive integers.
We regard the polynomial ring $S = K[x, y, z]$ as a ${\Bbb Z}$-graded ring
by $\deg(x) = a$, $\deg(y) = b$ and $\deg(z) = c$.

We denote by ${\Bbb P}_K(a,b,c)$ the weighted projective space $\proj S$.
Let
\[
\pi : X_K(a,b,c) \longrightarrow  {\Bbb P}_K(a,b,c)
\]
be the blow-up at the point corresponding to ${\frak p}_K(a,b,c)$.
We remark that ${\Bbb P}_K(a,b,c)$ is non-singular at this point
(e.g., Lemma~9 in \cite{Cu}).
If no confusion is possible, we denote 
${\frak p}_K(a,b,c)$ (resp.\ $X_K(a,b,c)$, ${\Bbb P}_K(a,b,c)$)
simply by ${\frak p}$ (resp.\ $X$, ${\Bbb P}$).
Let $E$ be the exceptional divisor of $\pi$.
Let $A$ be a Weil divisor on $X$ which satisfies 
${\mathcal O}_X(A) = \pi^*{\mathcal O}_{\Bbb P}(1)$.
Since $a$, $b$, $c$ are pairwise coprime,
we have ${\mathcal O}_X(nA) = \pi^*{\mathcal O}_{\Bbb P}(n)$ 
for any $n \in {\Bbb Z}$ (e.g., \cite{Mori}).
Then, 
\[
{\rm Cl}(X) = {\Bbb Z}A + {\Bbb Z}E \simeq {\Bbb Z}^2
\]
with the intersection pairing
\[
A^2 = \frac{1}{abc}, \ \ E^2 = -1, \ \ A.E = E.A = 0 .
\]

\begin{Definition}\label{Knegative curve}
A curve $C$ on $X_K(a,b,c)$ is called a {\em negative curve} on $X_K(a,b,c)$
if $C^2 < 0$ and $C \neq E$.

An irreducible homogeneous polynomial $\xi$ in $[{\frak p}_K(a,b,c)^{(r)}]_{d}$ is called
a {\em negative curve} in ${\frak p}_K(a,b,c)^{(r)}$ if $d/r < \sqrt{abc}$.
\end{Definition}

If a negative curve $C$ on $X_K(a,b,c)$ exists, then it is unique.
If a negative curve $\xi$ in  $[{\frak p}_K(a,b,c)^{(r)}]_{d}$ exists,
then $r$ and $d$ are uniquely determined, and $\xi$ 
is also unique up to multiplication
by an element in $K^\times$.
The negative curve $C$ on $X_K(a,b,c)$ is the proper transform
of $V_+(\xi)$.

\begin{Lemma}\label{Kproper}
Let $K$ be a field and $a$, $b$, $c$ be pairwise coprime positive integers.
We assume that  ${\frak p}_K(a, b, c)$ is minimally generated by
the three elements in (\ref{Kseisei}).

Then the curve $V_+(z^u - x^{s_3}y^{t_3})$  in ${\Bbb P}_K(a,b,c)$ 
is isomorphic to ${\Bbb P}^1_K$.
The proper transform $C$ (in $X$) of this curve  is also isomorphic to ${\Bbb P}^1_K$.
\end{Lemma}

\proof
First of all, we remark that $z^u - x^{s_3}y^{t_3}$ is an irreducible polynomial
by definition of $u$ (see  (\ref{Kstu})).
We put $v = x^{s_2}z^{u_2}/y^t$ and $w = x^{s_3}y^{t_3}/z^u$.

Since ${\frak p}_K(a, b, c)$ is generated by the three elements as in (\ref{Kseisei}),
we have 
\[
S[x^{-1}, y^{-1}, z^{-1}]_0 =  K[v^{\pm 1}, w^{\pm 1}] .
\]
Then, we have
\begin{equation}\label{NSy^{-1}}
S[y^{-1}]_0 = 
K \left[ v^\alpha w^\beta \ \left| \ \alpha, \beta \in {\Bbb Z}, \ \ 
\alpha \ge 0, \ \ 
 -\frac{s_2}{s_3}\alpha \le \beta \le \frac{u_2}{u}\alpha \right. \right] .
\end{equation}
Taking the degree $0$ component of 
\[
\frac{S[y^{-1}]}{(z^u - x^{s_3}y^{t_3}) S[y^{-1}]} \subset
\frac{S[x^{-1}, y^{-1}, z^{-1}]}{(w - 1) S[x^{-1}, y^{-1}, z^{-1}]} ,
\]
we obtain
\[
\frac{S[y^{-1}]_0}{(w - 1) K[v^{\pm 1}, w^{\pm 1}]  \cap S[y^{-1}]_0} \subset
\frac{K[v^{\pm 1}, w^{\pm 1}]}{(w - 1) K[v^{\pm 1}, w^{\pm 1}]} .
\]

Let $\phi: K[v^{\pm 1}, w^{\pm 1}] \rightarrow K[v^{\pm 1}]$ be the map
given by $\phi(w) = 1$.
The kernel of the map $\phi$ is $(w - 1) K[v^{\pm 1}, w^{\pm 1}]$.
By (\ref{NSy^{-1}}), we have $\phi(S[y^{-1}]_0) = K[v]$.
Hence, we have 
\[
\left[
\frac{S[y^{-1}]}{(z^u - x^{s_3}y^{t_3}) S[y^{-1}]}
\right]_0 \simeq K[v] .
\]
In the same way, we know that 
\[
\left[
\frac{S[x^{-1}]}{(z^u - x^{s_3}y^{t_3}) S[x^{-1}]}
\right]_0
\]
is also isomorphic to a polynomial ring over $K$ with one variable.
Hence, the curve $V_+(z^u - x^{s_3}y^{t_3})$  in ${\Bbb P}_K(a,b,c)$ 
is isomorphic to ${\Bbb P}^1_K$.

Since the map $C \rightarrow  V_+(z^u - x^{s_3}y^{t_3})$ is a finite birational map,
 $C$  is also isomorphic to ${\Bbb P}^1_K$.
\qed

\begin{Lemma}\label{Kckp}
Let $K$ be a field of prime characteristic $p$.
Let $a$, $b$, $c$ be pairwise coprime positive integers.
We assume the conditions (ii) and (iii) in Theorem~\ref{Kmain1}.

Then, there exist integers $q_1$ and $q_2$ such that
\[
H^i(X, {\mathcal O}_X(mA - nE)) = 0
\]
for $i > 0$, $n \ge q_1u$ and $m \ge (ab/u)n + q_2$. 
\end{Lemma}

\proof
First of all, remember that ${\mathcal O}_X(mA - nE)$ is invertible if $abc$ divides $m$ (e.g., Lemma~1.3 in \cite{Mori}).
By the condition (iii), $z^u - x^{s_3}y^{t_3}$ is the negative curve in ${\frak p}$.
Letting $C$ be the proper transform of $V_+(z^u - x^{s_3}y^{t_3})$,
$C$ is the negative curve that is linearly equivalent to $cuA-E$.
Therefore, $mA - nE$ is a nef Cartier divisor 
if $m \ge (ab/u)n \ge 0$ and $abc \mid m$.

Let $m$ and $n$ be integers such that $m \ge (ab/u)n \ge 0$.
Then, there exists a nef Cartier divisor $m_1A - n_1E$ such that
$0 \le m-m_1 < abc$ and $0 \le n-n_1 < cu$.

We use Fujita's vanishing theorem (Theorem~1.4.35 in \cite{La})
for finitely many coherent shaves 
\[
\left\{ {\mathcal O}_X(m_2A - n_2E) \mid  
0 \le m_2 < abc, \ 0 \le n_2 < cu \right\} .
\]
Then, there exists an ample Cartier divisor $q_2A + q_1(abA-uE)$ such that
\[
H^i(X, {\mathcal O}_X((m_2A - n_2E) + (m_1A - n_1E) + q_2A + q_1(abA-uE))) = 0
\]
for $i > 0$, $0 \le m_2 < abc$, $0 \le n_2 < cu$ and any nef Cartier divisor $m_1A - n_1E$.
Then, $q_1$ and $q_2$ satisfy the requirement in Lemma~\ref{Kckp}.
\qed

\begin{Lemma}\label{Kchp}
Let $K$ be a field of prime characteristic $p$.
Let $a$, $b$, $c$ be pairwise coprime positive integers.
We assume the conditions (ii) and (iii) in Theorem~\ref{Kmain1}.

Then, there exist $e > 0$ and $\eta \in [{\frak p}^{(p^{e}u)}]_{p^{e}ab}$ such that
$z^u - x^{s_3}y^{t_3}$ and $\eta$ satisfy Huneke's condition on ${\frak p}$, that is,
\[
\length{S}{S/(x, z^u -  x^{s_3}y^{t_3}, \eta)} = p^{e}u \cdot \length{S}{S/(x) + {\frak p}}
\]
holds. 
(The above integer $e$ depends on $a$, $b$, $c$ and $p$.)
\end{Lemma}

\proof 
Let $C$ be the proper transform of $V_+(z^u - x^{s_3}y^{t_3})$.
Then, $C$ is the negative curve on $X$ that is linearly equivalent to $ucA - E$. 

Consider the reflexive sheaf ${\mathcal O}_{\Bbb P}(ab)$.
Since $S_{ab}$ contains both $x^b$ and $y^a$, 
${\mathcal O}_{\Bbb P}(ab)$ is invertible away from the point $V_+(x,y)$.
Therefore ${\mathcal O}_X(abA)$ is invertible away from the point
$\pi^{-1}(V_+(x,y))$.
Since $C$ does not contain the point $\pi^{-1}(V_+(x,y))$, 
${\mathcal O}_X(abA - uE) \otimes {\mathcal O}_{nC}$ is 
an invertible sheaf on $nC$ for any $n > 0$. 

We choose integers $q_1$ and $q_2$ that satisfy Lemma~\ref{Kckp}.
By the condition (iii) in Theorem~\ref{Kmain1},
we obtain 
\[
uc < \frac{ab}{u} .
\]
Let $q$ be an integer which satisfies
\begin{equation}\label{Kdefofq}
qu \left( \frac{ab}{u} - uc \right) > q_2 .
\end{equation}

Consider the invertible sheaf ${\mathcal O}_X(abA - uE) \otimes {\mathcal O}_{C}$. 
Since $(abcA - cuE).C = 0$, the degree of ${\mathcal O}_X(abcA - cuE) \otimes {\mathcal O}_{C}$ is $0$.
Since $C$ is isomorphic to ${\Bbb P}^1_K$, 
${\mathcal O}_X(abcA - cuE) \otimes {\mathcal O}_{C} \simeq {\mathcal O}_{C}$.
Therefore, 
\begin{equation}\label{Korder}
{\mathcal O}_X(abA - uE) \otimes {\mathcal O}_{C} \simeq {\mathcal O}_{C}
\end{equation}
 since ${\rm Pic}({\Bbb P}^1_K) \simeq {\Bbb Z}$.

We have the exact sequences
\[
0 \longrightarrow {\mathcal O}_X(\ell C)/{\mathcal O}_X((\ell + 1) C)
\longrightarrow {\mathcal O}_{(\ell + 1) C} \longrightarrow 
{\mathcal O}_{\ell C} \longrightarrow 0 
\]
for $\ell = 1, 2, \ldots, qu - 1$.
They induce the exact sequences~\footnote{Suppose that $I$ is an ideal of a ring $A$ with
$I^2 = (0)$.  Then, consider the map $I \rightarrow A^\times$ defined by $a \mapsto 1+a$.
It induces the exact sequence $0 \rightarrow I \rightarrow A^\times \rightarrow (A/I)^\times \rightarrow 1$.}
\[
0 \longrightarrow {\mathcal O}_X(\ell C)/{\mathcal O}_X((\ell + 1) C)
\longrightarrow {\mathcal O}_{(\ell + 1) C}^\times \longrightarrow 
{\mathcal O}_{\ell C}^\times \longrightarrow 1 
\]
and 
\[
\begin{array}{ccccc}
H^1(X, {\mathcal O}_X(\ell C)/{\mathcal O}_X((\ell + 1) C)) & \longrightarrow &
H^1(X, {\mathcal O}_{(\ell + 1) C}^\times) & \longrightarrow & 
H^1(X, {\mathcal O}_{\ell C}^\times) \\
& & \parallel &  & \parallel \\
& & {\rm Pic}((\ell + 1) C) & \longrightarrow &  {\rm Pic}( \ell C)
\end{array}
\]
for $\ell = 1, 2, \ldots, qu - 1$.
Therefore we know that the order of  an element in the kernel of the map
${\rm Pic}((\ell + 1) C) \rightarrow  {\rm Pic}( \ell C)$ is $1$ or $p$.
Hence, the order of ${\mathcal O}_X(abA - uE) \otimes {\mathcal O}_{quC}$
(in ${\rm Pic}(qu C)$) is a power of $p$ by (\ref{Korder}).
Thus, there exists $e > 0$ such that 
\begin{equation}\label{Kp^e}
p^e > q + q_1
\end{equation}
and
\[
{\mathcal O}_X(p^e(abA - uE)) \otimes {\mathcal O}_{quC} \simeq {\mathcal O}_{quC} .
\]
Since 
\[
0 \longrightarrow {\mathcal O}_X(-qu C)
\longrightarrow {\mathcal O}_X \longrightarrow 
{\mathcal O}_{qu C} \longrightarrow 0 
\]
is exact, we have the following exact sequence:
\[
\begin{array}{l}
0 \rightarrow  {\mathcal O}_X(p^e(abA \mathop{-} uE)\mathop{-}qu C) \rightarrow 
{\mathcal O}_X(p^e(abA \mathop{-} uE))  \rightarrow 
{\mathcal O}_X(p^e(abA \mathop{-} uE)) \otimes {\mathcal O}_{qu C}  \rightarrow  0 \\
\hphantom{0 \longrightarrow  {\mathcal O}_X(p^e(abA - uE)-qu C) \longrightarrow 
{\mathcal O}_X(p^e(abA - uE))  \longrightarrow 
{\mathcal O}_X(p^e(abA - } \parallel \\
\hphantom{0 \longrightarrow  {\mathcal O}_X(p^e(abA - uE)-qu C) \longrightarrow 
{\mathcal O}_X(p^e(abA - uE))  \longrightarrow 
{\mathcal O}_X(p^e(abA - } {\mathcal O}_{qu C}
\end{array}
\]
Since $C \sim cuA-E$, we have
\[
p^e(abA - uE)-qu C 
 \sim  p^e(abA - uE)-qu (cuA-E)  
 =  (p^eab - qu^2c)A - (p^e - q)uE .
\]
By (\ref{Kp^e}), we have $p^e - q > q_1$.
By (\ref{Kdefofq}), we have
\[
p^eab - qu^2c > p^eab - qu(ab/u) + q_2 = (ab/u)(p^e - q)u + q_2 .
\]
Then, by Lemma~\ref{Kckp}, we have
\[
H^1(X, {\mathcal O}_X(p^e(abA - uE)-qu C)) = 0 .
\]
Then, we have
\[
\begin{array}{ccccc}
H^0(X, {\mathcal O}_X(p^e(abA - uE))) & \longrightarrow &
H^0(X, {\mathcal O}_X(p^e(abA - uE)) \otimes {\mathcal O}_{qu C} ) 
& \longrightarrow & 0 . \\
& & \parallel & & \\
& & H^0(X, {\mathcal O}_{qu C}) & & 
\end{array}
\]
Therefore, the natural map
\[
{\mathcal O}_X(p^e(abA - uE)) \longrightarrow 
{\mathcal O}_X(p^e(abA - uE))  \otimes {\mathcal O}_C \simeq {\mathcal O}_C
\]
induces the surjection
\[
H^0(X, {\mathcal O}_X(p^e(abA - uE))) \longrightarrow H^0(C, {\mathcal O}_C) = K .
\]
Thus, there exists an effective Weil divisor $D$ such that
$D \sim p^e(abA - uE)$ and the support of $D$ does not intersect with $C$.
Let $\eta$ be the equation of $\pi(D)$.
The degree of $\eta$ is $p^eab$.
Since $D \cap C = \emptyset$, $V_+(z^u - x^{s_3}y^{t_3}) \cap V_+(\eta) \subset V_+({\frak p})$ as a set.
Therefore, ${\frak p}$ is the only one minimal prime ideal of 
$(z^u - x^{s_3}y^{t_3}, \eta)$.
Hence $x$, $z^u - x^{s_3}y^{t_3}$, $\eta$ form a regular sequence of $S$.
We obtain
\[ 
\length{S}{S/(x, z^u -  x^{s_3}y^{t_3}, \eta )} = 
\length{S}{S/(x, z^u, y^{p^ea})} = p^eau = p^eu \cdot \length{S}{S/(x) + {\frak p}} ,
\]
which is the required equality.
\qed

\vspace{0.5em}
\noindent
{\it Proof of Theorem~\ref{Kmain1}.}\,
By Lemma~\ref{Kchp}, we know that 
\[
up^e \in \hc{{\frak p}_{\fp}(a,b,c)}{1}{z^u -  x^{s_3}y^{t_3}}.
\]
Then, by Theorem~\ref{N2l}, we know that 
$\srees{{\frak p}_{\qq}(a,b,c)}$ is Noetherian if and only if
\[
u \in \hc{{\frak p}_{\qq}(a,b,c)}{1}{z^u -  x^{s_3}y^{t_3}}.
\]
Thus,
We have completed the proof of Theorem~\ref{Kmain1}.
\qed

\begin{Remark}
\begin{rm}
Let $a$, $b$, $c$ be pairwise coprime positive integers.

Let $\xi \in [{\frak p}_K(a,b,c)^{(k)}]_d$ be a negative curve with $d/k < \sqrt{abc}$.
Then, $\srees{{\frak p}_{K}(a,b,c)}$ is Noetherian if and only if
$\hc{{\frak p}_{K}(a,b,c)}{k}{\xi} \neq \emptyset$.

Let $\xi$ be a homogeneous element in $S_\zz$.
Then, $\xi_\qq \in [{\frak p}_\qq(a,b,c)^{(k)}]_d$ is a negative curve with $d/k < \sqrt{abc}$ if and only if, for $p \gg 0$,
$\xi_\fp \in [{\frak p}_\fp(a,b,c)^{(k)}]_d$ is a negative curve with $d/k < \sqrt{abc}$.

Let $\xi$ be a homogeneous element in $S_\zz$.
Assume that  $\xi_\qq \in [{\frak p}_\qq(a,b,c)^{(k)}]_d$ is a negative curve with $d/k < \sqrt{abc}$, and 
$\srees{{\frak p}_\qq(a,b,c)}$ is Noetherian.
Then, $\hc{{\frak p}_{\qq}(a,b,c)}{k}{\xi_\qq} = \hc{{\frak p}_{\fp}(a,b,c)}{k}{\xi_\fp}$
for $p \gg 0$.

Let $\xi \in [{\frak p}_K(a,b,c)^{(k)}]_d$ be a negative curve with $d/k < \sqrt{abc}$.
Assume that $k = 1$ or $2$.
Suppose that
\[
\xi \equiv y^i \,\,{\rm mod}\,\, xS
\]
for some $i$.
Furthermore, assume that $\srees{{\frak p}_K(a,b,c)}$ is Noetherian.
Then, there exists a positive integer $m$ such that
\[
\hc{{\frak p}_K(a,b,c)}{k}{\xi} = \{ \ell m \mid \ell \in {\Bbb N} \}
\]
by Proposition~\ref{N2i} (1).
\end{rm}
\end{Remark}

\section{The condition EU}\label{Ksec3}

In this section, we introduce a sufficient condition (which is called as ``the condition EU'' below) 
for finite generation of $\srees{{\frak p}}$
under the assumption in Theorem~\ref{Kmain1}.
The condition EU was defined in Ebina~\cite{Ebina} and Uchisawa~\cite{Uchi}. 
We shall prove that, if $u \le 6$, the condition EU is a necessary and sufficient condition 
for finite generation of $\srees{{\frak p}}$ in Proposition~\ref{Kmain2}.

\vspace{2mm}

Let us remember the method introduced in Gonz\'ales-Karu~\cite{GK}.
Let $a$, $b$, $c$ be pairwise coprime positive integers and $K$ be a field.
Let $S = K[x, y, z]$ be a ${\Bbb Z}$-graded ring with $\deg(x) = a$, $\deg(y) = b$ and $\deg(z) = c$.
Suppose that the prime ideal ${\frak p}_K(a,b,c)$ is minimally generated by the three elements in (\ref{Kseisei}).

We put 
\[
v = x^{s_2}z^{u_2}/y^t, \ \ w = x^{s_3}y^{t_3}/z^u.
\]
Since ${\frak p}_K(a, b, c)$ is generated by the three elements in (\ref{Kseisei}),
we have 
\[
S[x^{-1}, y^{-1}, z^{-1}]_0 =  K[v^{\pm 1}, w^{\pm 1}] .
\]
Therefore, for each non-negative integer $e$, we have
\begin{equation}\label{Kvw}
S[x^{-1}, y^{-1}, z^{-1}]_{eab} =  y^{ea} \cdot K[v^{\pm 1}, w^{\pm 1}] .
\end{equation}

Let $\Delta_u$ be the domain (with boundary) surrounded by the following three lines
\begin{eqnarray*}
y & = & -(s_2/s_3)x \\
y & = & (u_2/u)x \\
y & = & (t/t_3)(x - u) + u_2 .
\end{eqnarray*}
Let $(0,0)$, $(u, u_2)$, $(\delta_1, \delta_2)$ be the vertices of  $\Delta_u$.
Here, $\delta_1$ and $\delta_2$ may not be integers.

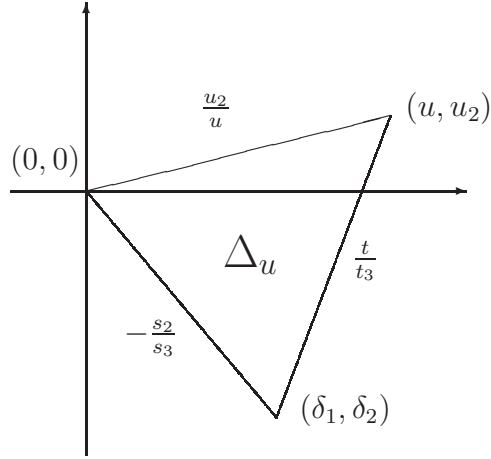
\[
{
\setlength\unitlength{1truecm}
  \begin{picture}(6,6.5)(-1,-4)
  \put(-1,0.3){$(0,0)$}
  \put(4.2,1){$(u,u_2)$}
  \put(2.8,-3){$(\delta_1,\delta_2)$}
  \put(1.8,-1){\mbox{\Large $\Delta_u$}}
  \put(-1,0){\vector(1,0){6}}
  \put(0,-3.5){\vector(0,1){6}}
  \put(0,0){\line(4,1){4}}
\qbezier (0,0) (1.25,-1.5) (2.5,-3)
\qbezier (4,1) (3.25,-1) (2.5,-3)
  \put(1.5,1){$\frac{u_2}{u}$}
  \put(3.5,-1){$\frac{t}{t_3}$}
  \put(0.5,-2){$-\frac{s_2}{s_3}$}
  \end{picture}
}
\]

For a non-negative integer $e$, we put
\[
e\Delta_u = \{ (e\alpha, e\beta) \mid (\alpha, \beta) \in \Delta_u \} .
\]
Then, it is easy to see that the euality~(\ref{Kvw}) induces
\[
S_{eab} =  y^{ea} \cdot \left(
\bigoplus_{(\alpha, \beta) \in e\Delta_u \cap {\Bbb Z}^2} K v^\alpha w^\beta
\right) .
\]
Since
\[
{\frak p}S[x^{-1}, y^{-1}, z^{-1}] = (v - 1, w - 1)S[x^{-1}, y^{-1}, z^{-1}] , 
\]
we have 
\[
{\frak p}^{(n)}S[x^{-1}, y^{-1}, z^{-1}] = (v - 1, w - 1)^nS[x^{-1}, y^{-1}, z^{-1}] 
\] 
and
\[
{\frak p}^{(n)} = (v - 1, w - 1)^nS[x^{-1}, y^{-1}, z^{-1}] \cap S
\]
for any $n > 0$.
Therefore,
\begin{equation}\label{KmethodGK}
\begin{array}{rcl}
[{\frak p}^{(n)}]_{eab} & = &  (v - 1, w - 1)^nS[x^{-1}, y^{-1}, z^{-1}] \cap S_{eab} \\
& = & {\displaystyle  
y^{ea} \left[
\left(
\bigoplus_{(\alpha, \beta) \in e\Delta_u \cap {\Bbb Z}^2} K v^\alpha w^\beta
\right) \bigcap
(v - 1, w - 1)^n K[v^{\pm 1}, w^{\pm 1}]
\right]  } .
\end{array}
\end{equation}

\begin{Remark}\label{Kch0}
\begin{rm}
Let $K$ be a field of characteristic $0$.
Let $\varphi(v,w)$ be an element in $K[v^{\pm 1}, w^{\pm 1}]$.
Then, $\varphi(v,w) \in (v - 1, w - 1)^n K[v^{\pm 1}, w^{\pm 1}]$ if and only if
\[
\frac{\partial^{k + \ell}\varphi}{\partial v^k \partial w^\ell}(1,1) = 0 
\]
for $k + \ell < n$.
\end{rm}
\end{Remark}

\begin{Remark}\label{Khunekegk}
\begin{rm}
Let $a$, $b$, $c$ be pairwise coprime positive integers.
Assume the conditions (ii) and (iii) in Theorem~\ref{Kmain1}.
Then, by Huneke's condition on ${\frak p}$, $\srees{{\frak p}}$ is Noetherian if and only if 
$[{\frak p}^{(eu)}]_{eab}$ contains an element whose coefficient of $y^{ea}$ is not $0$
for a sufficiently divisible $e$.

Assume the conditions (i),  (ii) and (iii) in Theorem~\ref{Kmain1}.
Then, by Theorem~\ref{Kmain1}, $\srees{{\frak p}}$ is Noetherian
if and only if $[{\frak p}^{(u)}]_{ab}$ contains an element whose coefficient of $y^a$ is not $0$.
By (\ref{KmethodGK}), it is equivalent to that 
\[
\left(
\bigoplus_{(\alpha, \beta) \in \Delta_u \cap {\Bbb Z}^2} K v^\alpha w^\beta
\right) \bigcap
(v - 1, w - 1)^u K[v^{\pm 1}, w^{\pm 1}]
\]
contains an element whose constant term is not $0$.
It is not so difficult to check whether it is satisfied or not using computers.
\end{rm}
\end{Remark}

Now, we introduce the condition EU which is defined by Ebina~\cite{Ebina} and Uchisawa~\cite{Uchi}.

\begin{Definition}[Ebina~\cite{Ebina}, Uchisawa~\cite{Uchi}]\label{Keu}
\begin{rm}
Let $a$, $b$, $c$ be pairwise coprime positive integers.
Suppose that the prime ideal ${\frak p}$ is minimally generated by the three elements in (\ref{Kseisei}).
For $i = 1, 2, \ldots, u$, we put
\[
\ell_i = ^\# \{ (\alpha, \beta) \in \Delta_u \cap {\Bbb Z}^2 \mid \alpha = i \} .
\]
We sort the sequence $\ell_1$, $\ell_2$, \ldots, $\ell_u$ into assending order
\[
\ell'_1 \le \ell'_2 \le \cdots \le \ell'_u .
\]

We say that the {\em condition EU} is satisfied for $(a,b,c)$ if
\[
\ell'_i \ge i
\]
for $i = 1, 2, \ldots, u$.
\end{rm}
\end{Definition}

\begin{Example}\label{Kexample}
\begin{rm}
\begin{itemize}
\item[(\RMN{1})]
Assume $(a,b,c) = (8,19,9)$.
Then, 
\[
{\frak p} = (x^7 - y^2z^2, y^3-x^6z, z^3 - xy) 
\]
and the conditions (ii) and (iii) in Theorem~\ref{Kmain1} is satisfied.
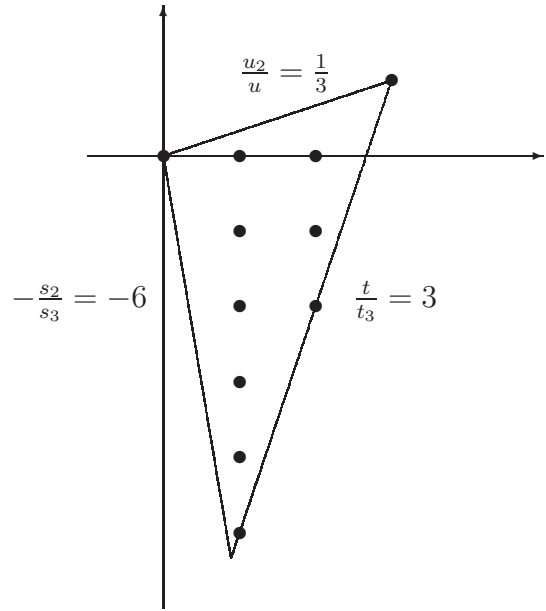
\[
{
\setlength\unitlength{1truecm}
  \begin{picture}(7,8)(-2,-6)
  \put(-1,0){\vector(1,0){6}}
  \put(0,-6){\vector(0,1){8}}
\qbezier (0,0) (1.5,0.5) (3,1)
\qbezier (0,0) (0.4444,-2.6666) (0.8889,-5.3333)
\qbezier (3,1) (1.0556,-4.8333) (0.8889,-5.3333)
\put(-0.1,-0.1){$\bullet$}
\put(0.9,-0.1){$\bullet$}
\put(0.9,-1.1){$\bullet$}
\put(0.9,-2.1){$\bullet$}
\put(0.9,-3.1){$\bullet$}
\put(0.9,-4.1){$\bullet$}
\put(0.9,-5.1){$\bullet$}
\put(1.9,-0.1){$\bullet$}
\put(1.9,-1.1){$\bullet$}
\put(1.9,-2.1){$\bullet$}
\put(2.9,0.9){$\bullet$}
  \put(1,1){$\frac{u_2}{u} = \frac{1}{3}$}
  \put(2.5,-2){$\frac{t}{t_3} = 3$}
  \put(-2,-2){$-\frac{s_2}{s_3} = -6$}
  \end{picture}
}
\]
Then, $u = 3$ and 
\[
\ell_1 = 6, \ \ \ell_2 = 3, \ \ \ell_3 = 1 .
\]
Therefore
\[
\ell'_1 = 1, \ \ \ell'_2 = 3, \ \ \ell'_3 = 6 .
\]
The condition EU is satisfied in this case.
\item[(\RMN{2})]
Assume $(a,b,c) = (25,29,72)$.
Then, 
\[
{\frak p} = (x^{11} - y^7z, y^{11}-x^7z^2,  z^3 - x^4y^4) 
\]
and the conditions (ii) and (iii) in Theorem~\ref{Kmain1} is satisfied.
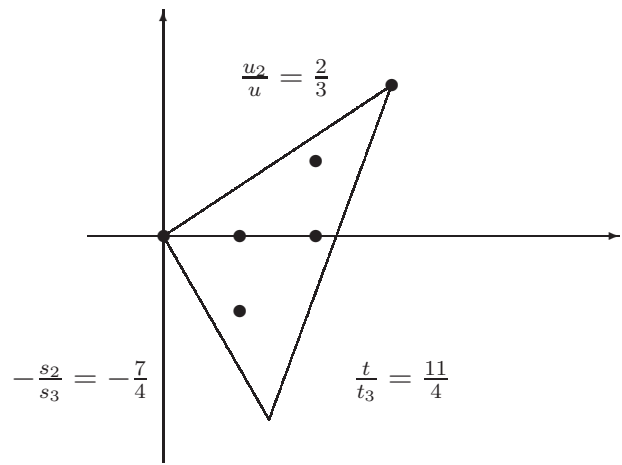
\[
{
\setlength\unitlength{1truecm}
  \begin{picture}(7,6)(-2,-3)
  \put(-1,0){\vector(1,0){7}}
  \put(0,-3){\vector(0,1){6}}
\qbezier (0,0) (1.5,1) (3,2)
\qbezier (0,0) (0.69445,-1.2153) (1.3889,-2.4306)
\qbezier (3,2) (2.19445,-0.2153) (1.3889,-2.4306)
\put(-0.1,-0.1){$\bullet$}
\put(0.9,-0.1){$\bullet$}
\put(0.9,-1.1){$\bullet$}
\put(1.9,0.9){$\bullet$}
\put(1.9,-0.1){$\bullet$}
\put(2.9,1.9){$\bullet$}
  \put(1,2){$\frac{u_2}{u} = \frac{2}{3}$}
  \put(2.5,-2){$\frac{t}{t_3} = \frac{11}{4}$}
  \put(-2,-2){$-\frac{s_2}{s_3} = -\frac{7}{4}$}
  \end{picture}
}
\]
Then, $u = 3$ and 
\[
\ell_1 = 2, \ \ \ell_2 = 2, \ \ \ell_3 = 1 .
\]
Therefore
\[
\ell'_1 = 1, \ \ \ell'_2 = 2, \ \ \ell'_3 = 2 .
\]
The condition EU is not satisfied in this case.
\item[(\RMN{3})]
Assume $(a,b,c) = (17,503,169)$.
Then, 
\[
{\frak p} = (x^{89} - y^2z^3, y^{3}-x^{49}z^4,  z^7 - x^{40}y) 
\]
and the conditions (ii) and (iii) in Theorem~\ref{Kmain1} is satisfied.
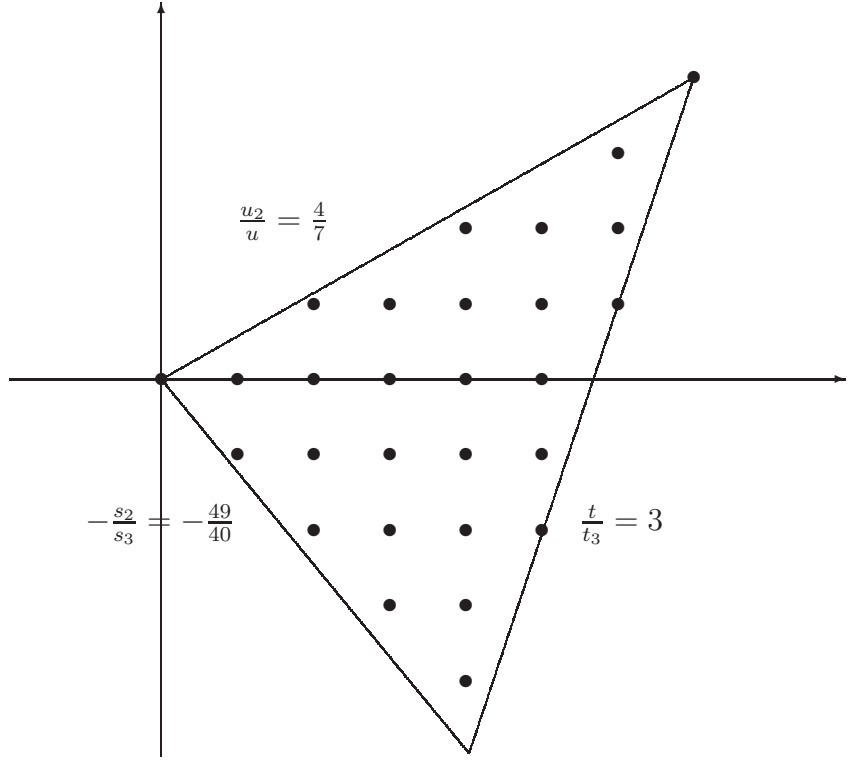
\[
{
\setlength\unitlength{1truecm}
  \begin{picture}(11,10)(-2,-5)
  \put(-2,0){\vector(1,0){11}}
  \put(0,-5){\vector(0,1){10}}
\qbezier (0,0) (3.5,2) (7,4)
\qbezier (0,0) (2.0238,-2.4792) (4.0476,-4.9583)
\qbezier (7,4) (5.5238,-0.47915)  (4.0476,-4.9583)
\put(-0.1,-0.1){$\bullet$}
\put(0.9,-0.1){$\bullet$}
\put(0.9,-1.1){$\bullet$}
\put(1.9,0.9){$\bullet$}
\put(1.9,-0.1){$\bullet$}
\put(1.9,-1.1){$\bullet$}
\put(1.9,-2.1){$\bullet$}
\put(2.9,0.9){$\bullet$}
\put(2.9,-0.1){$\bullet$}
\put(2.9,-1.1){$\bullet$}
\put(2.9,-2.1){$\bullet$}
\put(2.9,-3.1){$\bullet$}
\put(3.9,1.9){$\bullet$}
\put(3.9,0.9){$\bullet$}
\put(3.9,-0.1){$\bullet$}
\put(3.9,-1.1){$\bullet$}
\put(3.9,-2.1){$\bullet$}
\put(3.9,-3.1){$\bullet$}
\put(3.9,-4.1){$\bullet$}
\put(4.9,1.9){$\bullet$}
\put(4.9,0.9){$\bullet$}
\put(4.9,-0.1){$\bullet$}
\put(4.9,-1.1){$\bullet$}
\put(4.9,-2.1){$\bullet$}
\put(5.9,2.9){$\bullet$}
\put(5.9,1.9){$\bullet$}
\put(5.9,0.9){$\bullet$}
\put(6.9,3.9){$\bullet$}
  \put(1,2){$\frac{u_2}{u} = \frac{4}{7}$}
  \put(5.5,-2){$\frac{t}{t_3} = 3$}
  \put(-1,-2){$-\frac{s_2}{s_3} = -\frac{49}{40}$}
  \end{picture}
}
\]
Then, $u = 7$ and 
\[
\ell_1 = 2, \ \ \ell_2 = 4, \ \ \ell_3 = 5, \ \ \ell_4 = 7, \ \ \ell_5 = 5, \ \ \ell_6 = 3, \ \ \ell_7 = 1 .
\]
Therefore
\[
\ell'_1 = 1, \ \ \ell'_2 = 2, \ \ \ell'_3 = 3, \ \ \ell'_4 = 4, \ \ \ell'_5 = 5, \ \ \ell'_6 = 5, \ \ \ell'_7 = 7 .
\]
The condition EU is not satisfied in this case.
\end{itemize}
\end{rm}
\end{Example}

In order to show that the condition EU is a sufficient condition for finite generation of $\srees{{\frak p}}$
under some assumptions, 
we need the following lemma (\cite{Ebina}, \cite{Uchi}).
For the convenience of the reader, we give a proof of it here.

\begin{Lemma}[Ebina~\cite{Ebina}, Uchisawa~\cite{Uchi}]\label{Keulemma}
Let $K$ be a field of characteristic $0$ and $v$, $w$ be variables.

Let $u$ be a positive integer and $\alpha_1$, $\alpha_2$, \ldots, $\alpha_u$ be mutually distinct integers.
For $i = 1, 2, \ldots, u$, consider the integers $\beta_{i1}$, $\beta_{i2}$, \ldots, $\beta_{ii}$ satisfying 
\[
\beta_{i1} < \beta_{i2} < \cdots < \beta_{ii} .
\]
Put 
\[
T = \bigcup_{i = 1}^u \{ (\alpha_i, \beta_{i1}),  (\alpha_i, \beta_{i2}), \ldots, (\alpha_i, \beta_{ii}) \} \subset {\Bbb Z}^2 .
\]

Then, we have 
\[
\left(
\bigoplus_{(\alpha, \beta) \in T} Kv^\alpha w^\beta 
\right)
\bigcap 
(v-1, w-1)^uK[v^{\pm 1}, w^{\pm 1}] = 0 .
\]
\end{Lemma}

\proof
We shall prove it by induction on $u$.

If $u = 1$, then $^\#T = 1$.
It is easily verified in this case.

Assume $u \ge 2$.
Take 
\[
\varphi(v,w) \in 
\left(
\bigoplus_{(\alpha, \beta) \in T} Kv^\alpha w^\beta 
\right)
\bigcap 
(v-1, w-1)^uK[v^{\pm 1}, w^{\pm 1}] .
\]
Considering $v^{-\alpha_u}\varphi(v,w)$,
we may assume $\alpha_u = 0$.
Then, $\frac{\partial \varphi}{\partial v}$ satisfies all the assumptions with $u-1$.
Here, recall $\frac{\partial \varphi}{\partial v} \in (v-1, w-1)^{u-1}K[v^{\pm 1}, w^{\pm 1}]$ by Remark~\ref{Kch0}.
By induction, we obtain  $\frac{\partial \varphi}{\partial v} = 0$.
Therefore, we may suppose
\[
\varphi(v,w) = \sum_{j = 1}^u C_j w^{\beta_{uj}}
\]
where $C_1, \ldots, C_u \in K$.
Since $\varphi(v,w) \in (v-1, w-1)^uK[v^{\pm 1}, w^{\pm 1}]$, we have
\[
0 = \frac{\partial^k \varphi}{\partial w^k}(1,1) = 
\sum_{j = 1}^u C_j \beta_{uj} (\beta_{uj} - 1) \cdots (\beta_{uj} - k +1)
\]
for $k = 0, 1, \ldots, u-1$.
Then, we have
\[
\sum_{j = 1}^u C_j \beta_{uj}^k = 0
\]
for $k = 0, 1, \ldots, u-1$.
It is easy to see $C_1 = C_2 = \cdots = C_u = 0$.
\qed

\vspace{0.5em}
By this lemma, we can prove that the condition EU is a sufficient condition
for finite generation of $\srees{{\frak p}}$ under some assumptions.

\begin{Proposition}[Ebina~\cite{Ebina}, Uchisawa~\cite{Uchi}]\label{Keuprop}
Let $a$, $b$, $c$ be pairwise coprime positive integers.
Assume the conditions (i), (ii), (iii) in Theorem~\ref{Kmain1}.

If the condition EU is satisfied, then $\srees{{\frak p}_K(a,b,c)}$ is Noetherian.
\end{Proposition}

\proof
By the condition EU, we can choose a set $T$ as in Lemma~\ref{Keulemma}
which satisfies
\[
T \subset (\Delta_u - \{ (0,0) \} ) \cap {\Bbb Z}^2 .
\]
By (\ref{KmethodGK}), we obtain
\[
[{\frak p}^{(u)}]_{ab} =
y^{a} \left[
\left(
\bigoplus_{(\alpha, \beta) \in \Delta_u \cap {\Bbb Z}^2} K v^\alpha w^\beta
\right) \bigcap
(v - 1, w - 1)^u K[v^{\pm 1}, w^{\pm 1}] \right] .
\]
By this equality, we know that $[{\frak p}^{(u)}]_{ab}$ is defined
by $\frac{u(u+1)}{2}$ linear equations in 
$y^a \left( \oplus_{(\alpha, \beta) \in \Delta_u \cap {\Bbb Z}^2} K v^\alpha w^\beta \right)$.
We put $T' = T \cup \{ (0,0) \}$.
Recall that $T' \subset \Delta_u \cap {\Bbb Z}^2$ and $^\#T' = \frac{u(u+1)}{2}+1$.
Then, $[{\frak p}^{(u)}]_{ab}$ contains a non-zero element in the form
\[
\eta(x,y,z) = y^{a} \left(
\sum_{(\alpha, \beta) \in T'} C_{(\alpha, \beta)} v^\alpha w^\beta
\right) ,
\]
where $C_{(\alpha, \beta)} \in K$.

If $C_{(0,0)} = 0$, we have 
\[
\eta(x,y,z) \in y^a 
\left[
\left(
\bigoplus_{(\alpha, \beta) \in T} Kv^\alpha w^\beta 
\right)
\bigcap 
(v-1, w-1)^uK[v^{\pm 1}, w^{\pm 1}]
\right] = 0
\]
by Lemma~\ref{Keulemma}.
It is a contradiction.
Therefore, $C_{(0,0)} \neq 0$.
Then, 
\[
\length{S}{S/(x, z^u -  x^{s_3}y^{t_3}, \eta)} = ua = u \cdot \length{S}{S/(x) + {\frak p}}
\]
holds.
Hence, $\srees{{\frak p}}$ is Noetherian by Huneke's condition.
\qed

\vspace{0.5em}
The aim in the rest of this section is to prove the converse of Proposition~\ref{Keuprop} in the case $u \le 6$.

\begin{Definition}\label{Kn_im_i}
Let $a$, $b$, $c$ be pairwise coprime positive integers.
Assume the condition (ii) in Theorem~\ref{Kmain1}.

We define 
\[
n = ^\# \left\{ [-(s_2/s_3) , (u_2/u)] \cap {\Bbb Z} \right\} , \ \ 
m = ^\# \left\{ [(u_2/u), (t/t_3)] \cap {\Bbb Z} \right\} ,
\]
where $[ \ , \ ]$ is the closed interval.

We say that the {\em condition GK} is satisfied if one of the following two conditions 
is satisfied:
\begin{itemize}
\item[(\RMN{1})]
$^\# \left\{ (n-1)[(u_2/u), (t/t_3)] \cap {\Bbb Z} \right\} = n$ and 
$(u_2/u)n \not\in {\Bbb Z}$,
\item[(\RMN{2})]
$^\# \left\{ (m-1)[-(s_2/s_3) , (u_2/u)] \cap {\Bbb Z} \right\} = m$ and 
$(u_1/u)m \not\in {\Bbb Z}$.
\end{itemize}
\end{Definition}

We remark that the above condition (\RMN{1}) is satisfied for $a$, $b$, $c$
if and only if the above condition (\RMN{2}) is satisfied for $b$, $a$, $c$.

Let $a$, $b$, $c$ be pairwise coprime positive integers.
Assume the conditions (i),  (ii), (iii) in Theorem~\ref{Kmain1}.
If the condition GK is satisfied, then $\srees{{\frak p}_K(a,b,c)}$ is not Noetherian
by Theorem~1.2 in Gonz\'ales-Karu~\cite{GK}.

\begin{Proposition}\label{KGK}
Let $a$, $b$, $c$ be pairwise coprime positive integers.
Assume the conditions (i),  (ii), (iii) in Theorem~\ref{Kmain1}.

Then, the condition GK is satisfied if and only if 
one of the following five conditions is satisfied:
\begin{itemize}
\item[(GK1)]
$n = 1$,
\item[(GK2)]
$m = 1$,
\item[(GK3)]
$n = m = 2 < u$,
\item[(GK4)]
$3 \le n < u$, $m = 2$ and $^\# \left\{ (n-1)[(u_2/u), (t/t_3)] \cap {\Bbb Z} \right\} = n$,
\item[(GK5)]
$n = 2$, $3 \le m < u$ and $^\# \left\{ (m-1)[-(s_2/s_3) , (u_2/u)] \cap {\Bbb Z} \right\} = m$.
\end{itemize}
\end{Proposition}

\proof
Let $(\delta_1, \delta_2)$ be one of the vertices of $\Delta_u$ as in the beginning of this section.
First, we remark that, if $0 \le i < i+1 \le \delta_1$, then $\ell_{i+1} \ge \ell_i + (n-1)$.
In the same way,  if $\delta_1 \le i < i+1 \le u$, then $\ell_i \ge \ell_{i+1} + (m-1)$.
Thus, it is easy to see the following:
\begin{equation}\label{Kn,m>2}
\mbox{If $n\ge 3$ and $m\ge 3$,
then the condition EU is satisfied.}
\end{equation}
\begin{equation}\label{Kn=2mgeu}
\mbox{If $n = 2$ and $m\ge u$,
then the condition EU is satisfied.}
\end{equation}
\begin{equation}\label{Km=2ngeu}
\mbox{If $n\ge u$ and $m = 2$,
then the condition EU is satisfied.}
\end{equation}

Next, recall $s = s_2+s_3$, $t = t_1+t_3$ and $u = u_1+u_2$ by the condition (ii).
Then, we have
\begin{eqnarray*}
a & = & \length{S}{S/(x) + {\frak p}} = \length{S}{S/(x, y^t, z^u, y^{t_1}z^{u_1})} = tu - t_3u_2 , \\
b & = & \length{S}{S/(y) + {\frak p}} = \length{S}{S/(x^s, y, z^u, x^{s_2}z^{u_2})} = su - s_3u_1 , \\
c & = & \length{S}{S/(z) + {\frak p}} = \length{S}{S/(x^s, y^t, z, x^{s_3}y^{t_3})} = st - s_2t_1 .
\end{eqnarray*}
Since $a$ and $b$ are coprime, $u_1$, $u_2$ and $u$ are pairwise coprime.
Therefore, $(u_2/u)n \not\in {\Bbb Z}$ if and only if $n/u \not\in {\Bbb Z}$, and
$(u_3/u)m \not\in {\Bbb Z}$ if and only if $m/u \not\in {\Bbb Z}$.

It is easy to see that, if  the condition (GK$i$) is satisfied for some $i$, then  the condition GK is satisfied.

Conversely, assume that the condition GK is satisfied.
Then $\srees{{\frak p}_K(a,b,c)}$ is not Noetherian
by Theorem~1.2 in Gonz\'ales-Karu~\cite{GK}.
By Theorem~\ref{Keuprop} and (\ref{Kn,m>2}),
either $n < 3$ or $m < 3$ is satisfied.

If $n = 1$ (resp.\ $m = 1$), then (GK1) (resp.\ (GK2)) holds.

Suppose $n = 2$ and $m \ge 2$.
Since the condition EU is not satisfied,
we have $m < u$ by (\ref{Kn=2mgeu}).
If (\RMN{1}) of the condition GK is satisfied, then $n = m = 2 < u$, and therefore
(GK3) is satisfied.
If (\RMN{2}) of the condition GK is satisfied, then (GK3) or (GK5) is satisfied.

Suppose $n \ge 3$ and $m = 2$.
We know $n < u$ by (\ref{Km=2ngeu}).
Then (\RMN{2}) is not satisfied.
If (\RMN{1}) is satisfied, then (GK4) is satisfied.
\qed

\begin{Lemma}\label{Kle}
Let $a$, $b$, $c$ be pairwise coprime positive integers.
Assume the conditions (i),  (ii), (iii) in Theorem~\ref{Kmain1}.

\begin{itemize}
\item[1)]
Assume $n = 2$ and $3 \le m < u$.
If either $u_1 = 1$ or $u_2 = 1$ is satisfied,
then either the condition GK or EU is satisfied. 
\item[2)]
If $n = 2$ and $u > m \ge (u+1)/2$, then either the condition GK or EU is satisfied. 
\end{itemize}
\end{Lemma}

\proof
First of all, remark that the condition EU is satisfied for $a$, $b$, $c$ 
if and only if so for $b$, $a$, $c$.
Furthermore, the condition GK is satisfied for $a$, $b$, $c$ 
if and only if so for $b$, $a$, $c$.

First, we shall prove 1).
Assume $n = 2$, $3 \le m < u$ and $u_2 = 1$.
If (GK5) is not satisfied, then
\begin{equation*}
\begin{split}
\ell_1 =2, \ \ell_2 \ge 3, \ \ldots, \ \ell_{m-2} \ge m-&1, \ \ell_{m-1} \ge m+1, \ \ldots \ , \\
  &\ell_{2m-3} \ge 2m-1, \ \ell_{2m-2} \ge 2m+1, \ldots
\end{split}
\end{equation*}
and
\[
\ell_u = 1, \ \ell_{u-1} \ge m, \ \ell_{u-2} \ge 2m, \ \ell_{u-3} \ge 3m, \ldots 
\]
since $u_2 = 1$.
Thus, the condition EU is satisfied.

Next, assume $n = 2$, $3 \le m < u$ and $u_1 = 1$.
Considering $b$, $a$, $c$, we may assume $3 \le n < u$, $m=2$ and $u_2 = 1$.
Then, we have
\[
\ell_u = 1, \ \ell_{u-1} \ge 2, \ \ell_{u-2} \ge 4, \ \ell_{u-3} \ge 6, \ \ell_{u-4} \ge 8, 
\ldots
\]
and
\[
\ell_1 = n \ge 3, \ \ell_2 \ge 2n-1 \ge 5, \ \ell_3 \ge 3n-2 \ge 7, \ \ell_4 \ge 4n-3 \ge 9, \ldots .
\] 
In this case, the condition EU is always satisfied.

Here, we start to prove 2).
Assume that (GK5) is not satisfied.
Then, we have
\[
\ell_1 =2, \ \ell_2 \ge 3, \ \ldots, \ \ell_{m-2} \ge m-1, \ \ell_{m-1} \ge m+1, \ 
\ell_{m} \ge m+2, \ \ldots
\]
and
\[
\ell_u = 1, \ \ell_{u-1} \ge m, \ \ell_{u-2} \ge 2m-1 \ge u, \  \ldots .
\]
Thus, the condition EU is satisfied.
\qed

\begin{Proposition}\label{Kmain2}
Let $a$, $b$, $c$ be pairwise coprime positive integers.
Assume the conditions (i),  (ii), (iii) in Theorem~\ref{Kmain1}.

If $u \le 6$, then the condition EU is a necessary and sufficient condition
for finite generation of $\srees{{\frak p}_K(a,b,c)}$.

If $u \le 6$, then the condition GK is a necessary and sufficient condition
for infinite generation of $\srees{{\frak p}_K(a,b,c)}$.
\end{Proposition}

\proof
We shall prove that either the condition GK or EU is satisfied if $u \le 6$.

If $n = 1$ or $m = 1$, then  (GK1) or (GK2) is satisfied.

If $n \ge 3$ and $m \ge 3$, then EU is satisfied as in (\ref{Kn,m>2}).

If $u = n = m = 2$, then EU is satisfied.

If $u > n = m = 2$, then (GK3) is satisfied.

If $n = 2$ and $m \ge u$, then then EU is satisfied by (\ref{Kn=2mgeu}).

If $n \ge u$ and $m =2$, then then EU is satisfied by (\ref{Km=2ngeu}).

Now assume that $n = 2$ and $3 \le m < u$.
If $u > m \ge (u+1)/2$, then either the condition GK or EU is satisfied
by Lemma~\ref{Kle} 2). 
Assume $3 \le m < (u+1)/2$.
If $u \le 5$, then such $m$ does not exist.
Suppose $u = 6$ and $m = 3$.
Since $u$, $u_1$, $u_2$ are pairwise coprime, either $u_1$ or $u_2$ is $1$.
Then, by Lemma~\ref{Kle} 1), either the condition GK or EU is satisfied.

Assume $3 \le n < u$ and $m = 2$.
If $u \le 6$, we can prove that either the condition GK or EU is satisfied
in the same way as above.
\qed

\begin{Example}\label{Kexample2}
\begin{rm}
Let $K$ be a field of characteristic $0$.

Remember the three examples in Example~\ref{Kexample}.

Assume $(a,b,c) = (8,19,9)$.
In this case, $u = 3$, and the conditions (i), (ii) and (iii) in Theorem~\ref{Kmain1}
are satisfied.
Since the conditin EU is satisfied, $\srees{{\frak p}_K(a,b,c)}$ is Noetherian.

Assume $(a,b,c) = (25,29,72)$.
In this case, $u = 3$, and the conditions (i), (ii) and (iii) in Theorem~\ref{Kmain1}
are satisfied.
Since the conditin EU is not satisfied, $\srees{{\frak p}_K(a,b,c)}$ is not Noetherian.
Infinite generation of this ring was proved by Goto-Nishida-Watanabe~\cite{GNW}.

Assume $(a,b,c) = (17,503,169)$.
In this case, $u = 7$, and the conditions (i), (ii) and (iii) in Theorem~\ref{Kmain1}
are satisfied.
In this case, neither GK nor EU is satisfied.
Applying Theorem~\ref{Kmain1}, we know that  $\srees{{\frak p}_K(a,b,c)}$ is not Noetherian by calculation using computers (see Remark~\ref{Khunekegk}).
\end{rm}
\end{Example}

Let $a$, $b$, $c$ be pairwise coprime positive integers.
Assume the conditions (i),  (ii), (iii) in Theorem~\ref{Kmain1}.
We do not know any example of finitely generated $\srees{{\frak p}_K(a,b,c)}$
such that the condition EU is not satisfied.

\vspace{2em}

\section{An example having negative curve in the second symbolic power}\label{Nsec4}

Let $S = K[x, y, z]$, where $K$ is a field and $x, y, z$ are indeterminates.
We set $\gm = (x, y, z)S$ and $R = S_\gm$.
In this section, we first take positive integers
$s_2, s_3, t_1, t_3, u_1, u_2$ arbitrarily, and set
\[
f = x^s - y^{t_1}z^{u_1},\, g = y^t - x^{s_2}z^{u_2},\, h = z^u - x^{s_3}y^{t_3},
\]
where $s = s_2 + s_3, t = t_1 + t_3, u = u_1 + u_2$.
Moreover, we set
\[
a = t_3u_1 + t_1u,\, b = s_3u_2 + s_2u,\, c = s_2t_3 + s_3t.
\]
Let us regard $S$ as a $\zz$-graded ring by setting
\[
\deg x = a,\, \deg y = b,\, \deg z = c.
\]
Then, we can check directly that $f, g, h$ are all homogeneous.
We set $I = (f, g, h)S$ and $\ga = IR$.

\begin{Lemma}\label{N5a}
We have the following relations;
\begin{itemize}
\item[(1)]
$y^{t_3} f + z^{u_1}g + x^{s_2}h = 0$,
\item[(2)]
$z^{u_2}f + x^{s_3}g + y^{t_3}h = 0$.
\end{itemize}
\end{Lemma}

\noindent
{\it Proof.}\,
Since $f, g, h$ are the maximal minors of the matrix
\[
\left(
\begin{array}{lll}
y^{t_3} & z^{u_1} & x^{s_2} \\
z^{u_2} & x^{s_3} & y^{t_1}
\end{array}
\right),
\]
we get the relations stated above.
\qed

\begin{Lemma}\label{N5b}
The following assertions hold.
\begin{itemize}
\item[(1)]
$(x) + I$ is $\gm$-primary.
\item[(2)]
$\ass_S S / I = 
\assh_S S / I$.
\item[(3)]
$I_\gp$ is generated by $2$ elements for any $\gp \in \assh_S S / I$.
\item[(4)]
$\length{S}{S / (x) + I^{(n)}} = (n(n + 1) / 2)\cdot a$
for any $0 < n \in \zz$.
\item[(5)]
We have $I \subseteq \gp_K(a, b, c)$,
and the equality holds if $\mathrm{GCD}( a, b, c ) = 1$.
\end{itemize}
\end{Lemma}

\noindent
{\it Proof.}\,
(1)\,
This holds as $(x) + I$ contains $x, y^t$ and $z^u$.

(2)\,
We get this assertion by Hilbert-Burch's theorem.

(3)\,
Let us take any $\gp \in \ass_S S / I$.
Then, as $x \not\in \gp$, we have $h \in (f, g)S_\gp$ by Lemma~\ref{N5a} (1),
so $I_\gp = (f, g)S_\gp$.

(4)\,
Since $(x) + I = (x, y^t, y^{t_1}z^{u_1}, z^u)$,
we have $\mult{xR}{R / I} = \length{S}{S / (x) + I)} = a$.
Let us take any $0 < n \in \zz$.
Then,
\begin{equation*}
\begin{split}
\length{S}{S / (x) + I^{(n)}}  &= \length{R}{R / xR + \ga^{(n)}} = \mult{xR}{R / \ga^{(n)}} \\
     &\hspace{8ex}= \sum_{P \in \assh_R R / \ga}\length{R_P}{R_P / \ga_P^n}\cdot\mult{xR}{R / P}.
\end{split}
\end{equation*}
For any $P \in \assh_R R / \ga$, $\gr{\ga_P}$ is isomorphic to a polynomial ring
with $2$ variables over $R_P / \ga_P$, so
\begin{equation*}
\begin{split}
\length{R_P}{R_P / \ga_P^n} &= \sum_{i = 1}^n\length{R_P}{\ga_P^{i - 1} / \ga_P^i} 
 = \sum_{i = 1}^n i\cdot\length{R_P}{R_P / \ga_P}  \\
 &\hspace{8ex}= \frac{n(n + 1)}{2}\cdot\length{R_P}{R_P / \ga_P}.
\end{split}
\end{equation*}
Thus we get
\begin{equation*}
\begin{split}
\length{S}{S / (x) + I^{(n)}} &=
  \frac{n(n + 1)}{2}\!\!\!\!\sum_{P \in \assh_R R / \ga}
           \length{R_P}{R_P / \ga_P}\cdot\mult{xR}{R / P}  \\
  &\hspace{8ex}= \frac{n(n + 1)}{2}\cdot\mult{xR}{R / I}
        = \frac{n(n + 1)}{2}\cdot a.
\end{split}
\end{equation*}

(5)\,
We set $\gp = \gp_K(a, b, c)$.
Since $f, g, h$ are all homogeneous, we have $I \subseteq \gp$.
Hence, we have
\[
a = \length{S}{S / (x) + I} \geq \length{S}{S / (x) + \gp} = 
   \length{R}{R / xR + \gp R} = \mult{xR}{R / \gp R}.
\]
Now, we assume $\mathrm{GCD}( a, b, c ) = 1$.
Then, as is well known, we have $\mult{xR}{R / \gp R} = a$, so we see
\[
\length{S}{S / (x) + I)} = \length{S}{S / (x) + \gp},
\]
which means  $(x) + I = (x) + \gp$.
Then, we have
\[
\gp = \gp \cap ((x) + I) = (\gp \cap (x)) + I = x\gp + I,
\]
from which the equality $\gp = I$ follows.
\qed

\begin{Lemma}\label{N5c}
Suppose $s_2 > s_3$, $t_1 = t_3 = 1$ and $u_1 < u_2$.
Then, the following assertions hold.
\begin{itemize}
\item[(1)]
There exists $\xi \in I^{(2)}$ such that
\begin{itemize}
\item[(i)]
$x^{s_3}\xi = z^{u_2 - u_1}f^2 - gh$,
\item[(ii)]
$z^{u_1}\xi = x^{s_2 - s_3}h^2 - fg$, and
\item[(iii)]
$\xi \equiv y^3$ mod $(x)$.
\end{itemize}
\item[(2)]
$(x) + I^{(2)} =
(x, y^3, y^2z^{2u_1}, yz^{u + u_1}, z^{2u})$.
\end{itemize}
\end{Lemma}

\noindent
{\it Proof.}\,
(1)\,
>From the relations (1) and (2) of Lemma~\ref{N5a}, we get
\[
-yfh = z^{u_1}gh + x^{s_2}h^2
\hspace{2ex}
\mbox{and}
\hspace{2ex}
-yfh = z^{u_2}f^2 + x^{s_3}fg,
\]
respectively.
Hence, we have
\[
z^{u_1}gh + x^{s_2}h^2 = z^{u_2}f^2 + x^{s_3}fg,
\]
so we get
\[
x^{s_3}(x^{s_2 - s_3}h^2 - fg) = z^{u_1}(z^{u_2 - u_1}f^2 - gh).
\]
Since $x^{s_3}, z^{u_1}$ is a regular sequence on $S$,
there exists $\xi \in S$ such that
\[
x^{s_3}\xi = z^{u_2 - u_1}f^2 - gh
\hspace{2ex}
\mbox{and}
\hspace{2ex}
z^{u_1}\xi = x^{s_2 - s_3}h^2 - fg.
\]
The first equality implies $x^{s_3}\xi \in I^2$, so $\xi \in I^{(2)}$.
The second equality implies $z^{u_1}\xi \equiv -fg$ mod $(x)$,
so $z^{u_1}\xi \equiv yz^{u_1}\cdot y^2$ mod $(x)$
as $f \equiv -yz^{u_1}$ mod $(x)$ and $g \equiv y^2$ mod $(x)$.
Hence, we get $\xi \equiv y^3$ mod $(x)$ since $z^{u_1}$ is regular on $S / (x)$.

(2)\,
Since $(x) + I = (x, y^2, yz^{u_1}, z^u)$ and $u_1 < u_2$,
we have
\[
(x) + I^2 = (x, y^4, y^3z^{u_1}, y^2z^{2u_1}, yz^{u + u_1}, z^{2u}).
\]
We set $J = (\xi) + I^2 \subseteq I^{(2)}$.
Then, as
\[
(x) + J = (x, y^3, y^2z^{2u_1}, yz^{u + u_1}, z^{2u}),
\]
we have
\[
\length{S}{S / (x) + J} = 3(u + u_1) = 3a = \length{S}{S / (x) + I^{(2)}}
\]
by Lemma~\ref{N5b} (4).
Hence we get the required assertion.
\qed

\begin{Lemma}\label{N5d}
Suppose $s_2 > 2s_3$, $t_1 = t_3 = 1$ and $u_1 < u_2 < 2u_1$.
Then, the following assertions hold.
\begin{itemize}
\item[(1)]
There exists $\zeta \in I^{(3)}$ such that
\begin{itemize}
\item[(i)]
$x^{s_3}\zeta = f^3 + z^{2u_1 - u_2}h\xi$,
\item[(ii)]
$z^{u_2 - u_1}\zeta = f\xi + x^{s_2 - 2s_3}h^3$, and
\item[(iii)]
$\zeta \equiv -y^4z^{2u_1 - u_2}$ mod $(x)$.
\end{itemize}
\item[(2)]
$(x) + I^{(3)} = (x, y^5, y^4z^{2u_1 - u_2}, y^3z^u, y^2z^{u + 2u_1},
         yz^{2u + u_1}, z^{3u})$.
\item[(3)]
$S[ IT, I^{(2)}T^2, I^{(3)}T^3 ] \subsetneq \srees{I}$.
\end{itemize}
\end{Lemma}

\noindent
{\it Proof.}\,
(1)\,
>From the relations (i) and (ii) of Lemma~\ref{N5c}, we get
\[
fgh = z^{u_2 - u_1}f^3 - x^{s_3}f\xi
\hspace{2ex}\mbox{and}\hspace{2ex}
fgh = x^{s_2 - s_3}h^3 - z^{u_1}h\xi,
\]
respectively.
Hence, we have
\[
z^{u_2 - u_1}f^3 - x^{s_3}f\xi = x^{s_2 - s_3}h^3 - z^{u_1}h\xi,
\]
so we get
\[
z^{u_2 - u_1}(f^3 + z^{2u_1 - u_2}h\xi) = x^{s_3}(f\xi + x^{s_2 - 2s_3}h^3).
\]
Since $x^{s_3}, z^{u_2 - u_1}$ is a regular sequence on $S$,
there exists $\zeta \in S$ such that
\[
x^{s_3}\zeta = f^3 + z^{2u_1 - u_2}h\xi
\hspace{2ex}\mbox{and}\hspace{2ex}
z^{u_2 - u_1}\zeta = f\xi + x^{s_2 - 2s_3}h^3.
\]
The first equality implies $x^{s_3}\zeta \in II^{(2)} \subseteq I^{(3)}$,
so $\zeta \in I^{(3)}$ as $x^{s_3}$ is regular on $S / I^{(3)}$.
The second equality implies $z^{u_2 - u_1}\zeta \equiv f\xi$ mod $(x)$,
so $z^{u_2 - u_1}\zeta \equiv -yz^{u_1}\cdot y^3$ mod $(x)$ as
$f \equiv -yz^{u_1}$ mod $(x)$ and $\xi \equiv y^3$ mod $(x)$.
Hence, we get $\zeta \equiv -y^4z^{2u_1 - u_2}$ mod $(x)$
since $z^{u_2 - u_1}$ is regular on $S / (x)$.

(2)\,
By Lemma~\ref{N5c} (2), we have
\[
(x) + II^{(2)} =
(x, y^5, y^4z^{u_1}, y^3z^u, y^2z^{u + 2u_1}, yz^{2u + u_1}, z^{3u})
\]
as $u_1 < u_2 < 2u_1$.
We set $J = (\zeta) + II^{(2)} \subseteq I^{(3)}$.
Then, as
\[
(x) + J = 
(x, y^5, y^4z^{2u_1 - u_2}, y^3z^u, y^2z^{u + 2u_1}, yz^{2u + u_1}, z^{3u}),
\]
we have
\[
\length{S}{S / (x) + J} = 6(u + u_1) = 6a = \length{S}{S / (x) + I^{(3)}}
\]
by Lemma~\ref{N5b} (4).
Hence we get the required assertion.

(3)\,
It is enough to show
$I^{(2)}I^{(3)} \subsetneq I^{(5)}$.
(We have $(I^{(2)})^2 + II^{(3)} = I^{(4)}$,
which can be verified in the same way.)
By Lemma~\ref{N5c} (2) and Lemma~\ref{N5d} (2){, we have
\[
(x) + I^{(2)}I^{(3)} = (x) +
\left(\begin{array}{l}
y^8,\, y^7z^{2u_1 - u_2},\, y^6z^{\min\{u,\, 4u_1 - u_2\}},\, y^5z^{4u_1}, \\
y^4z^{4u_1 + u_2},\, y^3z^{3u},\, y^2z^{3u + 2u_1},\, yz^{4u + u_1},\, z^{5u}
\end{array}\right),
\]
so we get
\[
\length{S}{S / (x) + I^{(2)}I^{(3)}} = \min\{ 29u_1 + 16u_2, 32u_1 + 14u_2 \}.
\]
On the other hand, by Lemma~\ref{N5b} (4), we have
\[
\length{S}{S / (x) + I^{(5)}} = 15a = 30u_1 + 15u_2.
\]
Since $\min\{ 29u_1 + 16u_2, 32u_1 + 14u_2 \} - (30u_1 + 15u_2) =
\min\{ u_2 - u_1, 2u_1 - u_2 \} > 0$, we see
\[
\length{S}{S / (x) + I^{(2)}I^{(3)}} > \length{S}{S / (x) + I^{(5)}},
\]
which means $I^{(2)}I^{(3)} \subsetneq I^{(5)}$.
\qed

\vspace{0.5em}
In the rest of this section,
let us denote $S, I, \xi$ and $\zeta$ by
$S_K, I_K, \xi_K$ and $\zeta_K$, respectively,
in order to emphasize that the coefficient field is $K$.

\begin{Theorem}\label{N5e}
Let us choose any rational numbers $\alpha$ and $\beta$ such that
\[
1 < \alpha < \frac{5}{4}
\hspace{2ex}\mbox{and}\hspace{2ex}
2 < \beta < \frac{7}{3} - \frac{\alpha - 1}{2 - \alpha}.
\]
Moreover, we choose positive integers $s_2, s_3, u_1$ and $u_2$ such that
\[
\frac{s_2}{s_3} = \beta
\hspace{2ex}\mbox{and}\hspace{2ex}
\frac{u_2}{u_1} = \alpha.
\]
Then, setting $t_1 = t_3 = 1$,
we get the following assertions.
\begin{itemize}
\item[(1)]
$s_2 > 2s_3$ and $u_1 < u_2 < 2u_1$.
\item[(2)]
Let $0 \ll q \in \zz$.
We denote by $k$ the largest integer which is not bigger than $q / 3$.
Then we have
$(x^{s_3}, z^{2u_1 - u_2})^k \subseteq (x^{q(s_2 - 2s_3) + 1}, z^{q(u_2 - u_1)})$.
\item[(3)]
Let $p$ be any prime number. 
Then $3p^{e_p} \in \hc{I_\fp}{2}{\xi_\fp}$ for any $e_p \gg 0$.
\item[(4)]
$3 \not\in \hc{I_K}{2}{\xi_K}$ for any field $K$.
\item[(5)]
$\srees{I_\qq}$ is infinitely generated.
\end{itemize}
\end{Theorem}

\noindent
{\it Proof.}\,
(1)\,
These inequalities hold since
$s_2 / s_3 > 2$ and $1 < u_2 / u_1 < 2$.

(2)\,
Since $(x^{s_3}, z^{2u_1 - u_2})^k$ is generated my
\[ \{ x^{(k - i)s_3}z^{i(2u_1 - u_2)} \mid i = 0, 1, \dots, k \},
\]
it is enough to show that
\[
(k - i)s_3 \leq q(s_2 - 2s_3)
\hspace{2ex}\Longrightarrow\hspace{2ex}
i(2u_1 - u_2) \ge q(u_2 - u_1)
\]
holds for any $i = 0, 1, \dots, k$.
So, we suppose $(k - i)s_3 \leq q(s_2 - 2s_3)$, where $i = 0, 1, \dots, k$.
Then, dividing both sides of this inequality by $s_3$, we get
\[
k - i \leq q(\beta - 2).
\]
Here, we write $q = 3k + \ell$, where $\ell = 0, 1, 2$.
Then, as
\[
k - i \leq 3k(\beta - 2) + \ell(\beta - 2),
\]
we have
\[
i \geq k - 3k(\beta - 2) - \ell(\beta - 2) = k(7 - 3\beta) - \ell(\beta - 2).
\]
Hence, we get
\begin{eqnarray*}
i(2 - \alpha) - q(\alpha - 1) & \geq & \{k(7 - 3\beta) - \ell(\beta - 2)\}(2 - \alpha) - (3k + \ell)(\alpha - 1) \\
 & = & k\{(2 - \alpha)(7 - 3\beta) - 3(\alpha - 1)\} + m,
\end{eqnarray*}
where $m = -\ell\{(\beta - 2)(2 - \alpha) + (\alpha - 1)\}$.
Now, we notice that our assumption
$\alpha < 5 / 4$ and $\beta < 7 / 3 - (\alpha - 1) / (2 - \alpha)$ implies
\[
3\beta(2 - \alpha) < 7(2 - \alpha) - 3(\alpha - 1),
\]
so we see
\[
(2 - \alpha)(7 - 3\beta) - 3(\alpha - 1) > 0.
\]
Since $q \gg 0$, we have $k \gg 0$ too, so it follows that
\[
i(2 - \alpha) - q(\alpha - 1) > 0
\]
as $m$ is a bounded number.
Then, multiplying both sides of $i(2 - \alpha) > q(\alpha -1)$ by $u_1$, we get
\[
i(2u_1 - u_2) > q(u_2 - u_1).
\]

(3)\,
By Lemma~\ref{N5d} (ii), we have a relation
\[
z^{u_2 - u_1}\zeta_\fp - x^{s_2 - 2s_3}h^3 = f\xi_\fp
\]
in $S_\fp$.
We take $0 \ll e_p \in \zz$ and put $q = p^{e_p}$.
Then, we have
\[
z^{q(u_2 - u_1)}\zeta_\fp^q + (-1)^qx^{q(s_2 - 2s_3)}h^{3q} = f^q\xi_\fp^q.
\]
Here, we write $q = 3k + \ell$, where $\ell = 0, 1, 2$.
Then,
\[
f^q\xi_\fp^q = (f^3)^k\cdot f^\ell\xi_\fp^q \in (f^3)^k\cdot I^{(2q + \ell)}.
\]
The relation (i) of Lemma~\ref{N5d} means $f^3 \in (x^{s_3}, z^{2u_1 - u_2})I_\fp^{(3)}$.
Hence, by (2), we have
\[
(f^3)^k \in (x^{q(s_2 - 2s_3) + 1}, z^{q(u_2 - u_1)})I_\fp^{(3k)}.
\]
Thus we see
\[
f^q\xi_\fp^q \in (x^{q(s_2 - 2s_3) + 1}, z^{q(u_2 - u_1)})I_\fp^{(3q)},
\]
so there exist $\sigma_\fp, \tau_\fp \in I_\fp^{(3q)}$ such that
\[
z^{q(u_2 - u_1)}\zeta_\fp^q + (-1)^qx^{q(s_2 - 2s_3)}h^{3q} =
  x^{q(s_2 - 2s_3) + 1}\sigma_\fp + z^{q(u_2 - u_1)}\tau_\fp.
\]
Then, we have
\[
z^{q(u_2 - u_1)}\{ \zeta_\fp^q - \tau_\fp \} = 
    x^{q(s_2 - 2s_3)}\{ (-1)^{q + 1}h^{3q} + x\sigma_\fp \}.
\]
Since $x^{q(s_2 - 2s_3)}, z^{q(u_2 - u_1)}$ is a regular sequence on $S_\fp$,
there exists $\eta_\fp \in S_\fp$ such that
\[
x^{q(s_2 - 2s_3)}\eta_\fp = \zeta_\fp^q - \tau_\fp
\hspace{2ex}\mbox{and}\hspace{2ex}
z^{q(u_2 - u_1)}\eta_\fp = (-1)^{q + 1}h^{3q} + x\sigma_\fp.
\]
The first equality implies $x^{q(s_2 - 2s_3)}\eta_\fp \in I_\fp^{(3q)}$,
so we have $\eta_\fp \in I_\fp^{(3q)}$ as $x^{q(s_2 - s_3)}$ is regular on $S / I_\fp^{(3q)}$.
The second equality implies $z^{q(u_2 - u_1)}\eta_\fp \equiv (-1)^{q + 1}h^{3q}$ mod $xS_\fp$,
so $z^{q(u_2 - u_1)}\eta_\fp \equiv (-1)^{q + 1}z^{3qu}$ mod $xS_\fp$ as
$h \equiv z^u$ mod $xS_\fp$.
Hence, we get $\eta_\fp \equiv (-1)z^{2q(u + u_1)}$ mod $xS_\fp$ since
$z^{q(u_2 - u_1)}$ is regular on $S_\fp / xS_\fp$.
Then, we have
\begin{eqnarray*}
\length{S_\fp}{S_\fp / (x, \xi_\fp, \eta_\fp)} & = &
  \length{S_\fp}{S_\fp / (x, y^3, z^{2q(u + u_1)})} \\
 & = & 6q(u + u_1)  \\
 & = & 2\cdot 3q\cdot\length{S_\fp}{S_\fp / (x) + I_\fp},
\end{eqnarray*}
and hence $3q \in \hc{I_\fp}{2}{\xi_\fp}$.

(4)\,
If $3 \in \hc{I_K}{2}{\xi_K}$, by Proposition~\ref{N2i} (3), we have
\[
S_K[ I_KT, I_K^{(2)}T^2, I_K^{(3)}T^3 ] = \srees{I_K},
\]
which contradicts to Lemma~\ref{N5d} (3).

(5)\,
Let us notice that $\xi_\qq \in \zz[x, y, z]$.
Then, setting $k = 2$ and $r = 3$ in Theorem~\ref{N2l},
we see that $\srees{I_\qq}$ is not finitely generated.
\qed

\begin{Example}\label{N5f}
Let $\alpha = 6 / 5$ and $\beta = 49 / 24$,
which satisfy the assumptions on $\alpha$ and $\beta$ of Theorem~\ref{N5e}.
We set
\[
s_2 = 49m,\, s_3 = 24m,\, t_1 = t_3 = 1,\, u_1 = 5n
\hspace{1ex}\mbox{and}\hspace{1ex}
u_2 = 6n,
\]
where $m, n$ are coprime positive integers such that $m$ is odd
and $n$ is not a multiple of $97$.
Then, we have
\[
a = 16n,\,b = 683mn
\hspace{1ex}\mbox{and}\hspace{1ex}
c = 97m.
\]
Since $683$ and $97$ are prime numbers, we get
$\mathrm{GCD}( a, b, c ) = 1$.
Hence, by Lemma~\ref{N5b} (5) and Theorem~\ref{N5e} (5),
we see that $\srees{\gp_\qq(a, b, c)}$ is infinitely generated.
If $m = n = 1$,
then $a, b, c$ are pairwise coprime,
and one can check directly that $\xi_K$ is the negative curve
for any field $K$.
\end{Example}

\vspace{3mm}

\noindent
\begin{tabular}{l}
Kazuhiko Kurano \\
Department of Mathematics \\
Faculty of Science and Technology \\
Meiji University \\
Higashimita 1-1-1, Tama-ku \\
Kawasaki 214-8571, Japan \\
{\tt kurano@meiji.ac.jp} \\
{\tt http://www.isc.ac.jp/\~{}kurano}
\end{tabular}

\vspace{3mm}

\noindent
\begin{tabular}{l}
Koji Nishida \\
Institute of Management and Information Technologies \\
Chiba University \\
Yayoi-cho 1-33, Inage-ku \\
Chiba 263-8522, Japan \\
{\tt nishida@math.s.chiba-u.ac.jp} \\
\end{tabular}

\end{document}